\documentclass[12pt]{article}

\usepackage{amsthm,amsmath,amssymb,latexsym,color,cite}
\newtheorem{thm}{Theorem}[section]
\newtheorem{lem}[thm]{Lemma}
\newtheorem{cor}[thm]{Corollary}

\newtheorem{prop}[thm]{Proposition}
\newtheorem{defi}[thm]{Definition}
\theoremstyle{remark}
\newtheorem*{rmk}{Remark}

\numberwithin{equation}{section}

\renewcommand{\qed}{{\hfill\rule{4pt}{7pt}}\medskip}

\def\pf{\noindent {\it Proof.} }
\def\tree{{\rm tree}}

\def\im{{\rm imp}}
\def\d{{\rm gdes}}

\def\eld{{\rm eld}}
\def\deg{{\rm deg}}
\def\del{{\rm young}}

\def\chb{{\rm bdeg}}
\def\reld{{\rm reld}}
\def\der{\overline{\rm young }}

\begin{document}

\begin{center}
{\Large\bf A Generalization of the Ramanujan Polynomials and Plane Trees}
\end{center}
\vskip 2mm
\centerline{Victor J. W. Guo$^1$  and Jiang Zeng$^2$}

\begin{center}{\footnotesize
Institut Camille Jordan,
Universit\'e Claude Bernard (Lyon I),
F-69622, Villeurbanne Cedex, France} \\[5pt]
{\footnotesize
{\tt $^1$guo@math.univ-lyon1.fr,  http://math.univ-lyon1.fr/\textasciitilde{guo}}\\
{\tt $^2$zeng@math.univ-lyon1.fr,  http://math.univ-lyon1.fr/\textasciitilde{zeng}}
}
\end{center}

\vskip 0.5cm \noindent{\bf Abstract.}{\small Generalizing a sequence of Lambert,
Cayley and Ramanujan,
Chapoton has recently introduced a polynomial sequence $Q_n:=Q_n(x,y,z,t)$ defined by
$$
Q_1=1,\quad Q_{n+1}=[x+nz+(y+t)(n+y\partial_y)] Q_n.
$$
In this paper we prove Chapoton's conjecture on the duality
formula: $Q_n(x,y,z,t)=Q_n(x+nz+nt,y,-t,-z)$, and
answer his question about the combinatorial interpretation of $Q_n$.
Actually we give combinatorial
interpretations of these polynomials in terms of plane trees, half-mobile trees,
and forests of plane trees. Our approach also leads
to a general formula that unifies several known
results for enumerating trees and plane trees. }

\vskip 0.2cm
{\small \noindent{\bf Keywords}: Ramanujan polynomials, plane tree, half-mobile tree,
forest, general descent, elder vertex, improper edge

\vskip 0.2cm
\noindent{\bf MR Subject Classifications}: Primary 05A15; Secondary 05C05 }

\section{Introduction}
It is well-known that the Lambert function
$w=\sum_{n\geq 1}n^{n-1}y^n/n!$ (see \cite{Lam1})
is the solution to the functional equation
$we^{-w}=y$ with $w(0)=0$ and a formula of Cayley~\cite{Cayley}
says that $n^{n-1}$ counts the rooted labeled trees on $n$ vertices.
There are various generalizations of Cayley's formula. In particular,
 an interesting refinement of the sequence $n^{n-1}$ appeared
in Ramanujan's work
(see \cite{Berndt,DR,Zeng,CG}) and is related to Lambert's series as follows:
Differentiating $n$ times Lambert's function $w$ with respect to $y$
(see \cite[Lemma 6]{Zeng}) yields
\begin{equation}\label{eq:rama}
w^{(n)}=\frac{e^{nw}}{(1-w)^n}R_n\left(\frac{1}{1-w}\right),
\end{equation}
where $R_n$ is a polynomial of degree $n-1$ and satisfies the recurrence:
\begin{equation}\label{eq:lambert}
R_1=1,\quad R_{n+1}(y)=[n(1+y)+y^2\partial_y]R_n(y).
\end{equation}
It follows that
$R_2=1+y,\,
R_3=2+4y+3y^2$ and $
R_4=6+18y+25y^2+15y^3.
$
Clearly formula \eqref{eq:lambert} implies that $R_n(y)$
is a polynomial with nonnegative integral coefficients such that
$R_n(0)=(n-1)!$ and the leading coefficient
is $(2n-3)!!$. As the Lambert function is equivalent to $w^{(n)}(0)=n^{n-1}$,
we derive from \eqref{eq:rama} that
$R_n(1)=n^{n-1}$.

On the other hand, Shor~\cite{Shor} and Dumont-Ramamonjisoa~\cite{DR}
have proved independently that the coefficient of $y^k$ in $R_n(y)$
counts \emph{rooted labeled trees} on $n$ vertices with $k$ ``improper edges."
In a subsequent paper~\cite[Eq.~(26)]{Zeng} Zeng proved that if we set
$$
Z=\frac{e^{xw}-1}{x}=w+\frac{w^2}{2!}x+\frac{w^3}{3!}x^2+\cdots,
$$
then differentiating $n$ times $Z$ with respect to $y$ yields
\begin{equation}\label{eq:prefun}
Z^{(n)}=\frac{e^{(x+n)w}}{(1-w)^n}P_n\left(\frac{1}{1-w}, y\right),
\end{equation}
where $P_n:=P_n(x,y)$ is a polynomial defined by the recurrence relation:
\begin{equation*}
P_1=1,\quad P_{n+1}=[x+n+y(n+y\partial_y)]P_n.
\end{equation*}
The polynomials $P_n(x,y)$ are later
called the {\it Ramanujan polynomials}~\cite{CG,Ch1}.

Recently, Chapoton~\cite{Ch3} generalized $P_n$ to the polynomials $Q_n:=Q_n(x,y,z,t)$
as follows:
\begin{equation}\label{eq:recdef}
Q_1=1,\quad Q_{n+1}=[x+nz+(y+t)(n+y\partial_y)] Q_n.
\end{equation}
In the context of \emph{operads}, Chapoton~\cite[p.5]{Ch2}(see also \cite{Ch1}) conjectured that the coefficient of
$x^iy^{j-i}$ in $Q_n(1,x,0,y)$ is the \emph{dimension} of homogeneous component of degree
$(i,j)$ of the so-called Ramanujan operads $\text{Ram}(\{1,2,\ldots, n\})$.

Clearly these are \emph{homogeneous polynomials} in $x,y,z,t$
of degree $n-1$ and $z$ is just a homogeneous parameter. For example, we have $Q_2 =x+y+z+t$ and
$$
Q_3 =x^2+3xy+3xz+3xt+3y^2+4yz+5yt+2z^2+4zt+2t^2.
$$
We can easily derive explicit product formulae  of $Q_n$
for some special values. Indeed,
setting $t=-y$ in \eqref{eq:recdef} yields
\begin{align}
Q_n(x,y,z,-y)&=\prod_{k=1}^{n-1}(x+kz),\label{eq:special2}
\end{align}
while setting $y=0$ in \eqref{eq:recdef} leads to
\begin{equation}\label{eq:factor}
Q_n(x,0,z,t)=\prod_{k=1}^{n-1}(x+kz+kt).
\end{equation}
Further factorization formulae can be derived
from the following  \emph{duality formula}, which
was conjectured by Chapoton~\cite{Ch3}.
\begin{thm}\label{thm:conj}
For $n\geq 1$, there holds
\begin{equation}\label{eq:conj}
Q_n(x,y,z,t)=Q_n(x+nz+nt,y,-t,-z).
\end{equation}
\end{thm}

It follows from \eqref{eq:special2} and \eqref{eq:conj} that when $y=z$
the polynomials $Q_n$ factorize completely into linear factors:
\begin{equation}\label{eq:qnxt}
Q_n(x,z,z,t)=\prod_{k=1}^{n-1}(x+nz+kt).
\end{equation}
In particular, we obtain
$
Q_n(1,1,1,1)=n!C_n,
$
where $C_n=\frac{1}{n+1}{2n\choose n}$ is the $n$-th \emph{Catalan number}.
 This leads us to first look for
a combinatorial interpretation in the set of \emph{labeled plane trees}
on $n+1$ vertices rooted at 1, of which the cardinality is $n!C_n$
 (see \cite[p.~220]{Stanley99}).
To this end,
as $Q_n(x,y,z,t)=z^{n-1}Q_n(x/z,y/z,1,t/z)$,
 it is convenient to write
$Q_{n}(x,y,1,t)$ as follows:
\begin{equation}\label{eq:expansion}
Q_n(x,y,1,t)=\sum_{k=0}^{n-1} Q_{n,k}(x,t)y^k.
\end{equation}
Now identifying the
coefficients of $y^k$ in \eqref{eq:recdef} we obtain
$Q_{1,0}(x,t) =1$ and for $n\geq 2$:
\begin{equation}\label{rec-qnkxt}
Q_{n,k}(x,t)=[x+n-1+t(n+k-1)]Q_{n-1,k}(x,t)+(n+k-2)Q_{n-1,k-1}(x,t),
\end{equation}
 where
$Q_{n,k}(x,t)=0$ if $k\geq n$ or $k<0$.
The first values of $Q_{n,k}(x,t)$ are given in
Table~\ref{t:qnkxt}.
\begin{table}\label{t:qnkxt}
\caption{Values of $Q_{n,k}(x,t)$.}
{\scriptsize
\begin{center}
\begin{tabular}{|l|r|r|r|r|}\hline
 $ k\backslash n$ & $1$ & $2$ & $3$ & $4$ \\ \hline
 $0$ & $1$ & $x+1+t$ & $x^2+3x+2+(3x+4)t+2t^2$
 & $x^3+6x^2+11x+6+(6x^2+22x+18)t+(11x+18)t^2+6t^3$ \\ \hline
 $1$ & & $1$ & $3x+4+5t$
 & $6x^2+22x+18+(26x+43)t+26t^2$ \\ \hline
 $2$ & & & $3$ & $15x+25+35t$ \\ \hline
 $3$ & & & & $15$ \\ \hline
 $\sum_k$ & $1$ &$x+2+t$& $(x+3+t)(x+3+2t)$ & $(x+4+t)(x+4+2t)(x+4+3t)$
 \\ \hline
\end{tabular}
\end{center}
}
\end{table}

In the next section we shall prove that the
polynomials $Q_{n,k}(x,t)$ count plane trees of $k$ improper edges
with respect to two new statistics ``eld" and ``young."
It turns out that the polynomials $Q_n$ are
the counterpart of Ramanujan's polynomials $P_n$ (which count rooted trees)
for plane trees.

 Originally Chapoton asked for a combinatorial
interpretation of $Q_{n}$ in the model of \emph{half-mobile trees},
we shall answer his question in Section~3 by establishing
a bijection from plane trees to half-mobile trees.

In Section~4 we shall unify and generalize several
classical formulae for the enumeration of trees and plane trees. For instance,
in Theorem~\ref{thm:pplane} we prove that
\begin{align*}
\sum_{T}t^{\eld(T)}\prod_{i=1}^{n}x_i^{\del_T(i)}
=\prod_{k=0}^{n-2}(x_1+\cdots+x_n+kt),
\end{align*}
where $T$ ranges over all plane tree trees on $\{1,\ldots,n\}$ and
$\eld$ and $\del$ are mentioned as before.
In Section~5 we give another combinatorial interpretation of
$Q_{n,k}(x,t)$ in terms of forests of plane trees,
which extends a previous result of Shor~\cite{Shor}.
In Section~6 we give a short proof of Theorem~\ref{thm:conj}.
We end this paper with some open problems.

\section{Combinatorial interpretations in plane trees}
Let $\pi=a_1\cdots a_n$ be a permutation of a totally ordered set of $n$ elements.
Recall that an element $a_i$ is said to be a {\it right-to-left minimum} of $\pi$
 if $a_i<a_j$ for every $j>i$. For our purpose we need to
introduce a \emph{dual} statistic on permutations as follows.
The integer $i$ ($1\leq i\leq n-1$) is called a {\it general descent} of $\pi$,
if there exists a $j>i$ such that $a_j<a_i$. In other words, the
general descents of $\pi$ are positions that
do not correspond to right-to-left minima.
The number of general descents of $\pi$ is denoted by $\d(\pi)$.
For instance, if $\pi=3\,6\,1\,4\,5\,8\,7$, then the general descents of $\pi$ are
$1, 2$ and $6$, so $\d(\pi)=3$.
Let $\mathfrak{S}_n$ denote the set of all permutations of the set $[n]:=\{1,\ldots, n\}$.
It is easy to see that the following identity holds:
\begin{equation}\label{eq:general}
\sum_{\pi\in\mathfrak{S}_n}t^{\d(\pi)}=(1+t)(1+2t)\cdots(1+(n-1)t).
\end{equation}

Throughout this paper, unless indicated otherwise,
 all trees are {\it rooted labeled trees} on a linearly ordered vertex set.
Given two vertices $i$ and $j$ of a tree $T$ we say that $j$ is a {\it descendant} of $i$
if the path from the root to $j$ passes through $i$.
In particular, each vertex is a descendant of itself.
Let $\beta_{T}(i)$ be the smallest descendant of $i$.
Furthermore, if $j$ is a descendant of $i$ and is also connected to $i$ by an edge,
then we say that $j$ is a \emph{child} of $i$ and denote the corresponding edge by
$e=(i,j)$, and if $j'$ is another child of $i$, then we call $j'$ a {\it brother} of $j$.

A {\it plane tree} (or {\it ordered tree}) is a rooted tree in which
the children of each vertex are linearly ordered. From now on, by saying that
$v_1,\ldots, v_m$ are all the children of a vertex $v$ of a plane tree $T$ we mean
that $v_i$ is the $i$-th child of $v$, counting from left to right.
A vertex $j$ of a plane tree $T$
is called {\it elder} if $j$ has
a brother $k$ to its right such that $\beta_T(k)<\beta_T(j)$; otherwise we say that
$j$ is \emph{younger}.
Note that the rightmost child of any vertex is always younger.
For any vertex $v$ of a plane tree $T$,
let $\eld_T(v)$ be the number of elder children of $v$ in $T$.
Clearly, we have
$$
\eld_T(v)=\d(\beta_T(v_1)\cdots\beta_T(v_m)),
$$
where $v_1,\ldots,v_m$ are all the children of $v$.
Let $\eld(T)$ be the number of elder vertices of $T$.
Clearly, if $T$ is a plane tree on $n$
 vertices with $n\geq 2$, then $\eld(T)\leq n-2$.
\begin{defi}
Let $e=(i,j)$ be an edge of a tree $T$. We say that $e$ is a {\it proper} edge
or $j$ is a {\it proper child} of $i$,
if $j$ is an elder child of $i$ or $i<\beta_{T}(j)$.
Otherwise, we say that  $e$ is an {\it improper} edge and $j$ is an {\it improper child}
of $i$.
\end{defi}

For example,
for the plane tree $T$ in Figure~\ref{fig-good}
the set of elder vertices is  $\{3,8,9,11,12,13\}$ and the set of the improper edges is
$\{(3,14), \,(4,1),\, (6,5), \,(10,4), \,(14,2),\,(14,7)\}$.

\medskip

Given a vertex $v$ in a tree $T$, denote by $\deg(v)$ or $\deg_{T}(v)$
 the number of children of $v$, then the number of younger children of $v$ is given by
$$
\del(v)=\del_T(v)=\deg_T(v)-\eld_T(v).
$$

Denote by $\mathcal{P}_{n,k}$ (respectively $\mathcal{O}_{n,k}$ )
the set of plane trees (respectively plane trees with root $1$) on $[n]$ with $k$ improper edges.
Moreover, we may impose some conditions on the sets $\mathcal P_{n,k}$ and $\mathcal{O}_{n,k}$
to denote the subsets of plane trees that satisfy these conditions. For example,
$\mathcal P_{n,k} [\deg(n)=0]$ stands for the subset of
$\mathcal P_{n,k}$ subject to the condition $\deg(n)=0$.

\setlength{\unitlength}{1mm}
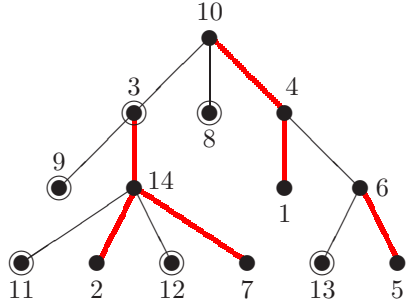
\begin{figure}[ht]
\begin{center}
{\footnotesize

\begin{picture}(50,37)(15,0)

\put(40,30){\line(0,-1){10}}
\linethickness{0.5mm}
\put(40,30){\line(-1,-1){10}}
\textcolor{red}{\qbezier(40,30)(45,25)(50,20)}
\put(30,20){\line(-1,-1){10}}
\put(30,20){\textcolor{red}{\line(0,-1){10}}}
\put(50,20){\textcolor{red}{\line(0,-1){10}}}
\put(50,20){\line(1,-1){10}}
\textcolor{red}{\qbezier(25,0)(27.5,5)(30,10)}
\put(30,10){\line(-3,-2){15}}
\put(30,10){\line(1,-2){5}}
\textcolor{red}{\qbezier(30,10)(37.5,5)(45,0)}
\put(60,10){\line(-1,-2){5}}
\textcolor{red}{\qbezier(60,10)(62.5,5)(65,0)}
\put(40,30){\circle*{2}} \put(40,30){\makebox(0,7)[c]{10}}
\put(30,20){\circle*{2}} \put(30,20){\makebox(0,7)[c]{3}}
\put(30,20){\circle{3}}
\put(40,20){\circle*{2}}
\put(40,20){\makebox(0,-7)[c]{8}}\put(40,20){\circle{3}}
\put(50,20){\circle*{2}} \put(50,20){\makebox(2,7)[c]{4}}
\put(20,10){\circle*{2}} \put(20,10){\makebox(0,7)[c]{9}}
\put(20,10){\circle{3}}
\put(30,10){\circle*{2}} \put(30,10){\makebox(7,2)[c]{14}}
\put(50,10){\circle*{2}} \put(50,10){\makebox(0,-7)[c]{1}}
\put(60,10){\circle*{2}} \put(60,10){\makebox(6,0)[c]{6}}
\put(15,0){\circle*{2}} \put(15,0){\makebox(0,-7)[c]{11}}
\put(15,0){\circle{3}}
\put(25,0){\circle*{2}} \put(25,0){\makebox(0,-7)[c]{2}}
\put(35,0){\circle*{2}} \put(35,0){\makebox(0,-7)[c]{12}}
\put(35,0){\circle{3}}
\put(45,0){\circle*{2}} \put(45,0){\makebox(0,-7)[c]{7}}
\put(55,0){\circle*{2}} \put(55,0){\makebox(0,-7)[c]{13}}
\put(55,0){\circle{3}}
\put(65,0){\circle*{2}} \put(65,0){\makebox(0,-7)[c]{5}}
\end{picture}
}

\end{center}

 \caption{A plane tree on $[14]$ with elder vertices circled and improper edges
thickened.\label{fig-good}}

\end{figure}

\begin{thm}\label{thm:qnkxt-1}
The polynomials $Q_{n,k}(x,t)$ have the following interpretation:
\begin{align}
Q_{n,k}(x,t)=\sum_{T\in \mathcal O_{n+1,k}}x^{\del_T(1)-1}t^{\eld(T)}.
\label{eq:qnkxt-1}
\end{align}
\end{thm}
\pf
Clearly, identity \eqref{eq:qnkxt-1}
is true for $n=1$. We shall prove by induction that the right-hand side of \eqref{eq:qnkxt-1}
satisfies the recurrence \eqref{rec-qnkxt} by distinguishing two cases according to
whether $n+1$ is a leaf or not.
\begin{itemize}
\item If $T\in\mathcal O_{n+1,k}$ with $\deg_T(n+1)=0$,
then deleting $n+1$ yields a plane tree $T'\in \mathcal O_{n,k}$.
Conversely, starting from any $T'\in \mathcal O_{n,k}$, we can recover $T$ by 
adding $n+1$ to $T'$ as a leaf
in $2n-1$ ways as follows. Pick up any vertex $v$ of $T'$ with the children being $a_1,\ldots,a_m$,
and then add $n+1$ as the $i$-th ($1\leq i\leq m+1$) child of $v$ to make the tree $T$. In other words,
the children of $v$ in $T$ become $a_1,\ldots,a_{i-1},n+1,a_{i+1},\ldots,a_m$.
Note that if $n+1$ is the rightmost child of $v$,
then $\eld(T)=\eld(T')$; otherwise, $\eld(T)=\eld(T')+1$.
Meanwhile, if $n+1$ is the rightmost child of $1$, then $\del_T(1)=\del_{T'}(1)+1$;
otherwise, $\del_T(1)=\del_{T'}(1)$.
 Since there are $n$ vertices in $T'$ and $\sum_{v\in [n]}\deg_{T'}(v)=n-1$, we obtain
\begin{align}
\sum_{T\in \mathcal O_{n+1,k}[\deg(n+1)=0]}&x^{\del_T(1)-1}t^{\eld(T)}\nonumber\\
&=[x+n-1+t(n-1)]\sum_{T\in \mathcal O_{n,k}}x^{\del_T(1)-1}t^{\eld(T)}.
 \label{eq:deg=0}
\end{align}
\item If $T\in\mathcal O_{n+1,k}$ with $\deg(n+1)>0$, suppose all the children of
$n+1$ are $a_1,\ldots,a_m$. Note that the edge $(n+1,a_{m})$ is always younger
and improper. We need to consider two cases:

If $m=1$ or $\beta_T(a_{m-1})<\beta_T(a_m)$, then replace $n+1$
by the child $a_m$ and contract the edge joining $n+1$
and $a_m$ such that the original children of $a_m$ are as the rightmost children
with their previous order unchanged. Thus, we obtain a plane tree $T'\in\mathcal O_{n,k-1}$.
Conversely, for such a $T'$, we can recover $T$ as follows. Pick any vertex
$v\neq 1$ of $T'$, and replace $v$ by $n+1$ and join $v$ to $n+1$ by an edge.
Suppose all the children of $v$ in $T'$ are $b_1,\ldots,b_p$, with
$b_{i_1},b_{i_2},\ldots,b_{i_r}$ ($i_1<i_2<\cdots< i_r$) being the improper children.
Then $\beta_{T'}(b_{i_j})<\beta_{T'}(b_s)$ for any $1\leq j\leq r$ and $s>i_j$.
It is easy to see that the children of $n+1$ in $T$
must be $b_1,b_2,\ldots,b_{i_j}$ and $v$, while the children of $v$ in $T$
must be $b_{i_j+1},b_{i_j+2},\ldots,b_p$, where $0\leq j\leq r$
and $i_0=0$. This means that there are $r+1$ possibilities to partition
the children of $v$. Since there are $k-1$ improper edges and $n$ vertices
(one of which is $1$) in $T'$, we see that there are total $n+k-2$
such corresponding plane trees $T$.
Moreover, we have $\del_{T}(1)=\del_{T'}(1)$ and
$\eld(T)=\eld(T')$. Hence the generating function for
the plane trees in $\mathcal O_{n+1,k}$ with $n+1$ having only one child or its second
child counting from right being younger is
\begin{align}
(n+k-2)\sum_{T\in \mathcal O_{n,k-1}}x^{\del_T(1)-1}t^{\eld(T)}.
\label{eq:bad}
\end{align}

If $m\geq 2$ and $\beta_T(a_{m-1})>\beta_T(a_m)$,
then $a_{m-1}$ is an elder child of $n+1$ and hence
$(n+1,a_{m-1})$ is a proper edge.
Replace $n+1$ by the child $a_{m-1}$ and contract the edge joining $n+1$
and $a_{m-1}$, such that the original children of $a_{m-1}$
are as the rightmost children in their previous order.
Thus, we obtain a plane tree $T'\in\mathcal O_{n,k}$.
Conversely, for such a $T'$, we can recover $T$ similarly as the first case.
Pick any vertex $v\neq 1$ of $T'$, and replace $v$ by $n+1$ and join $v$
to $n+1$ by an edge.
Suppose the children of $v$ in $T'$ are $b_1,\ldots,b_p$.
Assume all the improper children of $v$ in $T'$ are
$b_{i_1},b_{i_2},\ldots,b_{i_r}$.
It is easy to see that the only possible children of $n+1$ in $T$
are $b_1,b_2,\ldots b_{i_j-1},v,b_{i_j}$, while the children of $v$ in $T$ are
$b_{i_j+1},b_{i_j+2},\ldots,b_p$, where $1\leq j\leq r$.
Namely there are $r$ possibilities to partition the children of $v$.
Since $T'$ has $k$ improper edges, we see that there are total $k$ such
preimages $T$. Moreover, we have $\del_{T}(1)=\del_{T'}(1)$ and
$\eld(T)=\eld(T')+1$. Hence the generating function for the plane trees in
$\mathcal O_{n+1,k}$ with $n+1$ having at least two children
and its second child counting from right being elder is
\begin{align}
k\,t\sum_{T\in \mathcal O_{n,k}}x^{\del_T(1)-1}t^{\eld(T)}.
\label{eq:good}
\end{align}
Summing \eqref{eq:bad} and \eqref{eq:good} we obtain
\begin{align}
\sum_{T\in \mathcal O_{n+1,k}[\deg(n+1)>0]}x^{\del_T(1)-1}t^{\eld(T)}
&=k\,t\sum_{T\in \mathcal O_{n,k}}x^{\del_T(1)-1}t^{\eld(T)} \nonumber\\
&\quad+(n+k-2)\sum_{T\in \mathcal O_{n,k-1}}x^{\del_T(1)-1}t^{\eld(T)}.
 \label{eq:deg>0}
\end{align}

\end{itemize}
The proof then follows from summarizing identities \eqref{eq:deg=0} and
\eqref{eq:deg>0}.
\qed

The plane trees in $\mathcal O_{4,1}$ with their weights are
listed in Figure~\ref{fig-o41}.
\setlength{\unitlength}{0.8mm}
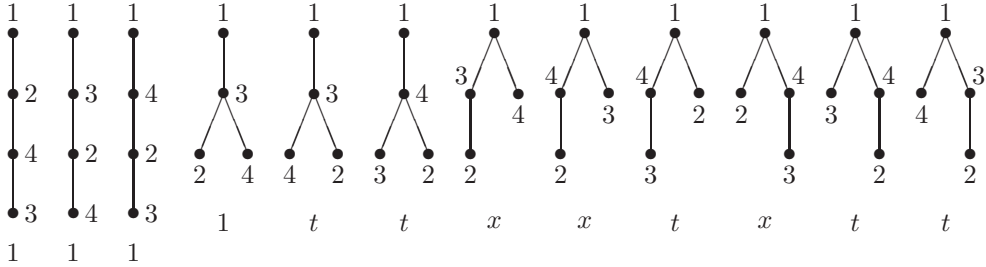
\begin{figure}[!hb]
\begin{center}
{\footnotesize
\begin{picture}(169,39)(0,21)
\put(5,55){\line(0,-1){30}}
\put(15,55){\line(0,-1){30}}
\put(25,55){\line(0,-1){30}}
\put(40,55){\line(0,-1){10}}
\put(40,45){\line(-2,-5){4}}
\put(40,45){\line(2,-5){4}}
\put(55,55){\line(0,-1){10}}
\put(55,45){\line(-2,-5){4}}
\put(55,45){\line(2,-5){4}}
\put(70,55){\line(0,-1){10}}
\put(70,45){\line(-2,-5){4}}
\put(70,45){\line(2,-5){4}}
\put(85,55){\line(-2,-5){4}}
\put(85,55){\line(2,-5){4}}
\put(81,45){\line(0,-1){10}}
\put(100,55){\line(-2,-5){4}}
\put(100,55){\line(2,-5){4}}
\put(96,45){\line(0,-1){10}}
\put(115,55){\line(-2,-5){4}}
\put(115,55){\line(2,-5){4}}
\put(111,45){\line(0,-1){10}}
\put(130,55){\line(-2,-5){4}}
\put(130,55){\line(2,-5){4}}
\put(134,45){\line(0,-1){10}}
\put(145,55){\line(-2,-5){4}}
\put(145,55){\line(2,-5){4}}
\put(149,45){\line(0,-1){10}}
\put(160,55){\line(-2,-5){4}}
\put(160,55){\line(2,-5){4}}
\put(164,45){\line(0,-1){10}}
\put(5,55){\circle*{1.6}}\put(5,55){\makebox(0,7)[c]{1}}
\put(5,45){\circle*{1.6}}\put(5,45){\makebox(6,0)[c]{2}}
\put(5,35){\circle*{1.6}}\put(5,35){\makebox(6,0)[c]{4}}
\put(5,25){\circle*{1.6}}\put(5,25){\makebox(6,0)[c]{3}}
\put(15,55){\circle*{1.6}}\put(15,55){\makebox(0,7)[c]{1}}
\put(15,45){\circle*{1.6}}\put(15,45){\makebox(6,0)[c]{3}}
\put(15,35){\circle*{1.6}}\put(15,35){\makebox(6,0)[c]{2}}
\put(15,25){\circle*{1.6}}\put(15,25){\makebox(6,0)[c]{4}}
\put(25,55){\circle*{1.6}}\put(25,55){\makebox(0,7)[c]{1}}
\put(25,45){\circle*{1.6}}\put(25,45){\makebox(6,0)[c]{4}}
\put(25,35){\circle*{1.6}}\put(25,35){\makebox(6,0)[c]{2}}
\put(25,25){\circle*{1.6}}\put(25,25){\makebox(6,0)[c]{3}}
\put(40,55){\circle*{1.6}}\put(40,55){\makebox(0,7)[c]{1}}
\put(40,45){\circle*{1.6}}\put(40,45){\makebox(6,0)[c]{3}}
\put(36,35){\circle*{1.6}}\put(36,35){\makebox(0,-7)[c]{2}}
\put(44,35){\circle*{1.6}}\put(44,35){\makebox(0,-7)[c]{4}}
\put(55,55){\circle*{1.6}}\put(55,55){\makebox(0,7)[c]{1}}
\put(55,45){\circle*{1.6}}\put(55,45){\makebox(6,0)[c]{3}}
\put(51,35){\circle*{1.6}}\put(51,35){\makebox(0,-7)[c]{4}}
\put(59,35){\circle*{1.6}}\put(59,35){\makebox(0,-7)[c]{2}}
\put(70,55){\circle*{1.6}}\put(70,55){\makebox(0,7)[c]{1}}
\put(70,45){\circle*{1.6}}\put(70,45){\makebox(6,0)[c]{4}}
\put(66,35){\circle*{1.6}}\put(66,35){\makebox(0,-7)[c]{3}}
\put(74,35){\circle*{1.6}}\put(74,35){\makebox(0,-7)[c]{2}}
\put(85,55){\circle*{1.6}}\put(85,55){\makebox(0,7)[c]{1}}
\put(81,45){\circle*{1.6}}\put(81,45){\makebox(-3,6)[c]{3}}
\put(89,45){\circle*{1.6}}\put(89,45){\makebox(0,-7)[c]{4}}
\put(81,35){\circle*{1.6}}\put(81,35){\makebox(0,-7)[c]{2}}
\put(100,55){\circle*{1.6}}\put(100,55){\makebox(0,7)[c]{1}}
\put(96,45){\circle*{1.6}}\put(96,45){\makebox(-3,6)[c]{4}}
\put(104,45){\circle*{1.6}}\put(104,45){\makebox(0,-7)[c]{3}}
\put(96,35){\circle*{1.6}}\put(96,35){\makebox(0,-7)[c]{2}}
\put(115,55){\circle*{1.6}}\put(115,55){\makebox(0,7)[c]{1}}
\put(111,45){\circle*{1.6}}\put(111,45){\makebox(-3,6)[c]{4}}
\put(119,45){\circle*{1.6}}\put(119,45){\makebox(0,-7)[c]{2}}
\put(111,35){\circle*{1.6}}\put(111,35){\makebox(0,-7)[c]{3}}
\put(130,55){\circle*{1.6}}\put(130,55){\makebox(0,7)[c]{1}}
\put(126,45){\circle*{1.6}}\put(126,45){\makebox(0,-7)[c]{2}}
\put(134,45){\circle*{1.6}}\put(134,45){\makebox(3,6)[c]{4}}
\put(134,35){\circle*{1.6}}\put(134,35){\makebox(0,-7)[c]{3}}
\put(145,55){\circle*{1.6}}\put(145,55){\makebox(0,7)[c]{1}}
\put(141,45){\circle*{1.6}}\put(141,45){\makebox(0,-7)[c]{3}}
\put(149,45){\circle*{1.6}}\put(149,45){\makebox(3,6)[c]{4}}
\put(149,35){\circle*{1.6}}\put(149,35){\makebox(0,-7)[c]{2}}
\put(160,55){\circle*{1.6}}\put(160,55){\makebox(0,7)[c]{1}}
\put(156,45){\circle*{1.6}}\put(156,45){\makebox(0,-7)[c]{4}}
\put(164,45){\circle*{1.6}}\put(164,45){\makebox(3,6)[c]{3}}
\put(164,35){\circle*{1.6}}\put(164,35){\makebox(0,-7)[c]{2}}
\put(5,20){\makebox(0,-3)[c]{$1$}}
\put(15,20){\makebox(0,-3)[c]{$1$}}
\put(25,20){\makebox(0,-3)[c]{$1$}}
\put(40,25){\makebox(0,-3)[c]{$1$}}
\put(55,25){\makebox(0,-3)[c]{$t$}}
\put(70,25){\makebox(0,-3)[c]{$t$}}
\put(85,25){\makebox(0,-3)[c]{$x$}}
\put(100,25){\makebox(0,-3)[c]{$x$}}
\put(115,25){\makebox(0,-3)[c]{$t$}}
\put(130,25){\makebox(0,-3)[c]{$x$}}
\put(145,25){\makebox(0,-3)[c]{$t$}}
\put(160,25){\makebox(0,-3)[c]{$t$}}
\end{picture}
}

\end{center}
\caption{The polynomial $Q_{3,1}(x,t)=3x+4+5t$ as a weight function of 
$\mathcal O_{4,1}$.\label{fig-o41}}
\end{figure}

\begin{table}
\caption{Values of $Q_{n,k}(x-t-1,t)$.}
{\scriptsize
\begin{center}
\begin{tabular}{|l|r|r|r|r|}\hline
 $ k\backslash n$ & $1$ & $2$ & $3$ & $4$ \\ \hline
 $0$ & $1$ & $x$ & $x^2+x+xt$
 & $x^3+3x^2+2x+(3x^2+4x)t+2xt^2$ \\ \hline
 $1$ & & $1$ & $3x+1+2t$
 & $6x^2+10x+2+(14x+7)t+6t^2$ \\ \hline
 $2$ & & & $3$ & $15x+10+20t$ \\ \hline
 $3$ & & & & $15$ \\ \hline
 $\sum_k$ & $1$ &$x+1$ & $(x+2)(x+2+t)$ & $(x+3)(x+3+t)(x+3+2t)$
 \\ \hline
\end{tabular}
\end{center}
}
\end{table}
\begin{thm}\label{thm:x-t-1}
The polynomials $Q_{n,k}(x-t-1,t)$ have the following interpretation:
\begin{align}
Q_{n,k}(x-t-1,t)=\sum_{T\in \mathcal P_{n,k}}x^{\del_T(1)}t^{\eld(T)}.
\label{eq:x-t-1}
\end{align}
\end{thm}
\pf
By \eqref{rec-qnkxt}, we see that
\begin{equation}\label{rec-qnkxt2} 
Q_{n,k}(x-t-1,t)=[x+n-2+t(n+k-2)]Q_{n-1,k}(x,t)+(n+k-2)Q_{n-1,k-1}(x,t).
\end{equation}
Similarly to the proof of Theorem \ref{thm:qnkxt-1}, we show that the
right-hand side of \eqref{eq:x-t-1} satisfies the recurrence \eqref{rec-qnkxt2}.
More precisely, we can prove that for $n\geq 1$,
\begin{align}
\sum_{T\in \mathcal P_{n,k}[\deg(n)=0]}x^{\del_T(1)}t^{\eld(T)}
&=[x+n-2+t(n-2)]\sum_{T\in \mathcal P_{n-1,k}}x^{\del_T(1)}t^{\eld(T)},
 \label{eq2:deg=0} \\
\sum_{T\in \mathcal P_{n,k}[\deg(n)>0]}x^{\del_T(1)}t^{\eld(T)}
&=k\,t\sum_{T\in \mathcal P_{n-1,k}}x^{\del_T(1)}t^{\eld(T)} \nonumber\\
&\quad+(n+k-2)\sum_{T\in \mathcal P_{n-1,k-1}}x^{\del_T(1)}t^{\eld(T)}.
 \label{eq2:deg>0}
\end{align}
The proof of \eqref{eq2:deg=0} and \eqref{eq2:deg>0} is exactly the same as that
of \eqref{eq:deg=0} and \eqref{eq:deg>0} and is omitted here. We only mention
that ``Pick any vertex $v\neq 1$ of $T'$" needs to be changed to
``Pick any vertex $v$ of $T'$," and if $n$ is the root of $T$ then we
take the vertex replacing $n$ as the root of $T'$. \qed

The plane trees in $\mathcal P_{3,1}$ with their weights are
listed in Figure~\ref{fig-p31}.

\setlength{\unitlength}{0.8mm}
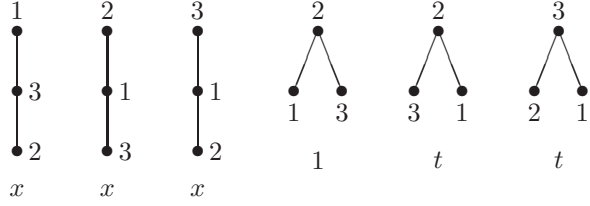
\begin{figure}[!h]
\begin{center}
{\footnotesize
\begin{picture}(104,31)(0,33)
\put(5,55){\line(0,-1){20}}
\put(20,55){\line(0,-1){20}}
\put(35,55){\line(0,-1){20}}
\put(55,55){\line(-2,-5){4}}
\put(55,55){\line(2,-5){4}}
\put(75,55){\line(-2,-5){4}}
\put(75,55){\line(2,-5){4}}
\put(95,55){\line(-2,-5){4}}
\put(95,55){\line(2,-5){4}}
\put(5,55){\circle*{1.6}}\put(5,55){\makebox(0,7)[c]{1}}
\put(5,45){\circle*{1.6}}\put(5,45){\makebox(6,0)[c]{3}}
\put(5,35){\circle*{1.6}}\put(5,35){\makebox(6,0)[c]{2}}
\put(20,55){\circle*{1.6}}\put(20,55){\makebox(0,7)[c]{2}}
\put(20,45){\circle*{1.6}}\put(20,45){\makebox(6,0)[c]{1}}
\put(20,35){\circle*{1.6}}\put(20,35){\makebox(6,0)[c]{3}}
\put(35,55){\circle*{1.6}}\put(35,55){\makebox(0,7)[c]{3}}
\put(35,45){\circle*{1.6}}\put(35,45){\makebox(6,0)[c]{1}}
\put(35,35){\circle*{1.6}}\put(35,35){\makebox(6,0)[c]{2}}
\put(55,55){\circle*{1.6}}\put(55,55){\makebox(0,7)[c]{2}}
\put(51,45){\circle*{1.6}}\put(51,45){\makebox(0,-7)[c]{1}}
\put(59,45){\circle*{1.6}}\put(59,45){\makebox(0,-7)[c]{3}}
\put(75,55){\circle*{1.6}}\put(75,55){\makebox(0,7)[c]{2}}
\put(71,45){\circle*{1.6}}\put(71,45){\makebox(0,-7)[c]{3}}
\put(79,45){\circle*{1.6}}\put(79,45){\makebox(0,-7)[c]{1}}
\put(95,55){\circle*{1.6}}\put(95,55){\makebox(0,7)[c]{3}}
\put(91,45){\circle*{1.6}}\put(91,45){\makebox(0,-7)[c]{2}}
\put(99,45){\circle*{1.6}}\put(99,45){\makebox(0,-7)[c]{1}}
\put(5,30){\makebox(0,-3)[c]{$x$}}
\put(20,30){\makebox(0,-3)[c]{$x$}}
\put(35,30){\makebox(0,-3)[c]{$x$}}
\put(55,35){\makebox(0,-3)[c]{$1$}}
\put(75,35){\makebox(0,-3)[c]{$t$}}
\put(95,35){\makebox(0,-3)[c]{$t$}}
\end{picture}
}
\end{center}
\caption{The polynomial $Q_{3,1}(x-t-1,t)=3x+1+2t$ as a weight function of
$\mathcal P_{3,1}$.\label{fig-p31}}
\end{figure}

It is worthwhile to point out that there is a simpler variant of
Theorems~\ref{thm:qnkxt-1} and \ref{thm:x-t-1}.
A vertex $j$ of a plane tree $T$ is called {\it really elder} if $j$ has
a brother $k$ to its right such that $k<j$. Let $\reld_T(v)$ be the number
of really elder children of $v$. Namely, $\reld_T(v)=\d(a_1\cdots a_m)$, where
$a_1,\ldots,a_m$ are all the children of $v$.
Assume $e=(i,j)$ is an edge of a tree $T$,
we say that $e$ is a {\it really proper} edge, if $j$ is a really elder child of
$i$ or $i<\beta_{T}(j)$. Otherwise, we call $e$ a {\it really improper} edge.
Let
$$\der_T(v)=\deg_T(v)-\reld_T(v),$$
and let $\overline{\mathcal P}_{n,k}$ (respectively $\overline{\mathcal O}_{n,k}$)
denote the set of plane trees (respectively plane trees with root $1$) on $[n]$ with
$k$ really improper edges.

\begin{cor}
There holds
\begin{align*}
Q_{n,k}(x,t)=\sum_{T\in \overline{\mathcal O}_{n+1,k}}x^{\der_T(1)-1}t^{\reld(T)}
=\sum_{T\in \overline{\mathcal P}_{n,k}}(x+t+1)^{\der_T(1)-1}t^{\reld(T)}.
\end{align*}
\end{cor}

\pf It suffices to construct a bijection $\phi$ from $\mathcal P_{n,k}$
(respectively $\mathcal O_{n+1,k}$)
to $\overline{\mathcal P}_{n,k}$ (respectively $\overline{\mathcal O}_{n+1,k}$).
Starting from a plane tree $T\in \mathcal P_{n,k}$ (respectively $T\in \mathcal O_{n+1,k}$),
we define $\phi(T)$ as the plane tree obtained from $T$ as follows.
For any vertex $v$ of $T$ with subtrees $T_1,\ldots, T_m$
rooted at $v_1,\ldots,v_m$, respectively, we reorder
these subtrees as
$T_{\sigma(1)},\ldots,T_{\sigma(m)}$ such that $v_{\sigma(i)}<v_{\sigma(j)}$
if and only if $\beta_T(v_{i})<\beta_T(v_{j})$, where $\sigma\in \mathfrak S_m$.
See Figure~\ref{fig:real}. \qed

\setlength{\unitlength}{1mm}
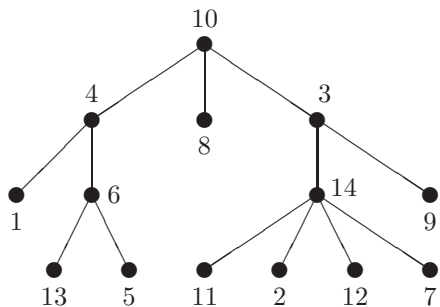
\begin{figure}[ht]
\begin{center}
{\footnotesize

\begin{picture}(50,37)(15,0)
\put(40,30){\line(-3,-2){15}}
\put(40,30){\line(0,-1){10}}
\put(40,30){\line(3,-2){15}}
\put(55,20){\line(3,-2){15}}
\put(55,20){\line(0,-1){10}}
\put(25,20){\line(-1,-1){10}}
\put(25,20){\line(0,-1){10}}
\put(55,10){\line(-3,-2){15}}
\put(55,10){\line(-1,-2){5}}
\put(55,10){\line(1,-2){5}}
\put(55,10){\line(3,-2){15}}
\put(25,10){\line(-1,-2){5}}
\put(25,10){\line(1,-2){5}}
\put(40,30){\circle*{2}} \put(40,30){\makebox(0,7)[c]{10}}
\put(25,20){\circle*{2}} \put(25,20){\makebox(0,7)[c]{4}}
\put(40,20){\circle*{2}} \put(40,20){\makebox(0,-7)[c]{8}}
\put(55,20){\circle*{2}} \put(55,20){\makebox(2,7)[c]{3}}
\put(70,10){\circle*{2}} \put(70,10){\makebox(0,-7)[c]{9}}
\put(55,10){\circle*{2}} \put(55,10){\makebox(7,2)[c]{14}}
\put(15,10){\circle*{2}} \put(15,10){\makebox(0,-7)[c]{1}}
\put(25,10){\circle*{2}} \put(25,10){\makebox(6,0)[c]{6}}
\put(40,0){\circle*{2}} \put(40,0){\makebox(0,-7)[c]{11}}
\put(50,0){\circle*{2}} \put(50,0){\makebox(0,-7)[c]{2}}
\put(60,0){\circle*{2}} \put(60,0){\makebox(0,-7)[c]{12}}
\put(70,0){\circle*{2}} \put(70,0){\makebox(0,-7)[c]{7}}
\put(20,0){\circle*{2}} \put(20,0){\makebox(0,-7)[c]{13}}
\put(30,0){\circle*{2}} \put(30,0){\makebox(0,-7)[c]{5}}
\end{picture}
}

\end{center}

\caption{The plane tree in $\overline{\mathcal P}_{14,6}$ corresponding
to that in Figure~\ref{fig-good}. \label{fig:real}}

\end{figure}

To end this section we make a connection to increasing (plane) trees.
A (plane) tree on $[n]$ is called {\em increasing} if any path from
the root to another vertex forms an increasing sequence.
Clearly an increasing plane tree has no improper edges and vice versa.
Combining \eqref{eq:factor}, \eqref{eq:expansion}
 and Theorem~\ref{thm:qnkxt-1} we get
 the following result.
\begin{prop}
For every $n\geq 1$, we have
\begin{align*}
\sum_{T\in \mathcal P_{n,0}}x^{\del_T(1)}t^{\eld(T)}
=\prod_{k=0}^{n-2}(x+k+kt).
\end{align*}
In particular,
the number of increasing trees on $[n]$ is $(n-1)!$ and the number of
increasing plane trees on $[n]$ is $(2n-3)!!$.
\end{prop}

\section{Combinatorial interpretations in half-mobile trees}
The notion of half-mobile trees was introduced by Chapoton~\cite{Ch2}.
A {\it half-mobile tree} on $[n]$ is defined to be a rooted tree
with two kinds vertices, called \emph{labeled} and \emph{unlabeled}
(or white and black, respectively) vertices satisfying the following conditions:
\begin{itemize}
\item The labeled vertices are in bijection with $[n]$;
\item Each unlabeled vertex has at least two children and all of them are
labeled;
\item There is a fixed cyclic order on the children of any unlabeled vertex;
\item The children of each labeled vertex are not ordered.
\end{itemize}
Let $T$ be a half-mobile tree. For any black vertex $x$ of $T$, we
define $\beta_T(x)$ to be the smallest vertex among
all the white descendants of $x$. From now on, we
assume that the rightmost child $x$ of a black vertex $v$ has the smallest
$\beta_T(x)$.
\begin{defi}
An edge $e=(u,v)$ of a half-mobile tree $T$
is called {\it improper} if $u$ and $v$ are labeled and $u>\beta_T(v)$,
or $u$ is unlabeled  and   $v$ is its rightmost child,
moreover $u$ has a (labeled) father  greater than $\beta_T(v)$.
\end{defi}

A {\it forest of half-mobile trees} on $[n]$ is a graph of which the
connected components are half-mobile trees and the white vertex set is
$[n]$. Denote by ${\mathcal H}_n$  the set of forests of half-mobile trees on $[n]$.
For any $F\in {\mathcal H}_n$, let $\im(F)$ be the number of improper edges
of $F$, and $\tree(F)$ the number of half-mobile trees of $F$.
Finally define the \emph{black degree} of $F$, denoted by $\chb(F)$, to be
the total degree of black vertices
minus the number of black vertices. Namely,
$$
\chb(F)=\sum_{v}(\deg_F(v)-1),
$$
where the sum is over all black vertices $v$ of $F$.

We first recall a fundamental transformation $\psi$ on
$\mathfrak{S}_n$. We identify each permutation $\pi\in \mathfrak{S}_n$
with the sequence $\pi(1)\pi(2)\ldots \pi(n)$.
 Since a variant of this bijection can be found in
\cite[Proposition 1.3.1]{Stanley97}, we only give an informal
 description of $\psi$ as follows:
\begin{itemize}
\item[(a)] Factorize the permutation $\pi$ into a product of disjoint cycles.
\item[(b)] Order the cycles of $\pi$ in increasing order of their minima.
\item[(c)] Write the minimum at last within each cycle
and erase the parentheses of cycles.
\end{itemize}

For instance, if the factorization of $\pi\in\mathfrak{S}_8$ is $(241)(73)(5)(86)$,
then $\psi(\pi)=2\,4\,1\,7\,3\,5\,8\,6$.

\begin{lem}\label{lemt}
The mapping $\psi$  is a bijection on $\mathfrak{S}_n$
 such that the number of cycles of $\pi$ is equal to the number
of right-to-left minima of $\psi(\pi)$.
\end{lem}

Now we are ready to construct our bijection from the plane trees to forests of half-mobile trees.
\begin{prop}\label{prop:theta}
There is a bijection $\theta: \mathcal O_{n+1}\longrightarrow {\mathcal H}_n$ such that
for any $T\in \mathcal O_{n+1}$
$$
\del_T(1)=\tree(\theta(T)),\quad
\eld(T)=\chb(\theta(T)), \quad\im(T)=\im(\theta(T)).
$$
\end{prop}

\pf Let $T\in\mathcal O_{n+1}$. Pick any vertex $v$ of $T$ with
$\deg(v)>0$. Suppose its children are $a_1,\ldots, a_m$.
Consider the permutation $\pi_v=b_1\cdots b_m$,
where $b_j=\beta_T(a_j)$ for $1\leq j\leq m$.
 Let $b_{i_1}, b_{i_2}, \ldots, b_{i_r}=b_m$
be all the right-to-left minima of $\pi_v$.
For any $k$ $(1\leq k\leq r)$, if $i_k-i_{k-1}>1$ ($i_0=0$), add a black vertex
$\rule{2mm}{2mm}$ to be a new child of $v$ and then move the subtrees rooted at
$a_{i_{k-1}+1},\ldots,a_{i_k}$ to be subtrees of this black vertex $\rule{2mm}{2mm}$ with
the cyclic order $(a_{i_{k-1}+1},\ldots,a_{i_k})$. In other words,
the cycle
$(a_{i_{k-1}+1},\ldots,a_{i_k})$ is a cycle of length greater than $1$
in the permutation $\psi^{-1}(\pi_v)$, obtained by applying the inverse mapping of
$\psi$ in Lemma~\ref{lemt}. Since all the elements except the rightmost one in
$(a_{i_{k-1}+1},\ldots,a_{i_k})$ are elder vertices in $T$,
applying the above procedure to every vertex  $v$  with $\deg(v)>0$,
we obtain a half-mobile tree $T'$ on $[n+1]$
with root $1$. By Lemma~\ref{lemt},
it is easy to see that the mapping $T\mapsto T'$ is reversible.
Furthermore,
 $$
\del_T(1)=\deg_{T'}(1),\quad
\eld(T)=\chb(T'), \quad\im(T)=\im(T').
$$
Deleting the root $1$, and shifting the label $i$ to $i-1$
($2\leq i\leq n+1$) in $T'$, we obtain a forest $\theta(T)$ of half-mobile trees on $[n]$.
This completes the proof.
\qed

As an illustration of the bijection $\theta$, we apply $\theta$ to a tree $T$ in
 Figure~\ref{fig:half0} with  $n=14$.
 For example, the root 1 of $T$ has three children $4,\, 8,\, 3$
 and $\pi_1=2\,8\,3$.
 Note that $\del_T(1)=2$, $\eld(T)=5$ and $\im(T)=4$.

\setlength{\unitlength}{1mm}
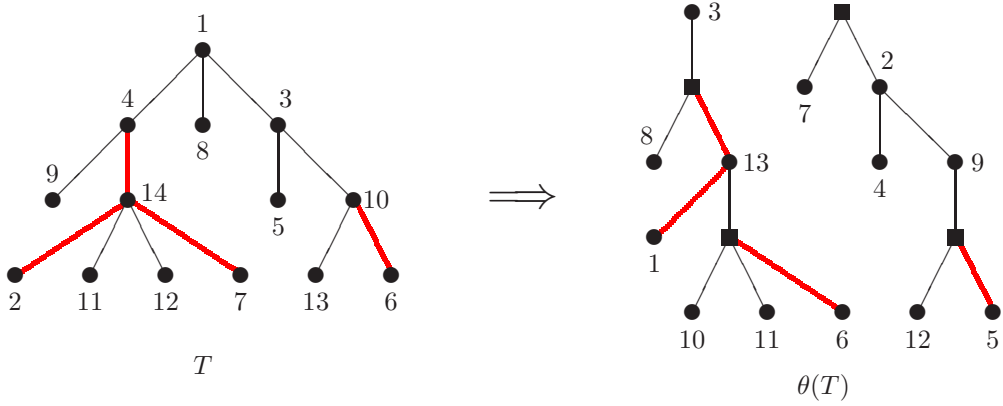
\begin{figure}[ht]
\begin{center}
{\footnotesize

\begin{picture}(130,49)(15,-10)

\put(40,30){\line(0,-1){10}}
\put(50,20){\line(0,-1){10}}
\linethickness{0.5mm}
\put(40,30){\line(-1,-1){10}}
\put(40,30){\line(1,-1){10}}
\put(30,20){\line(-1,-1){10}}
\put(30,20){\textcolor{red}{\line(0,-1){10}}}
\put(50,20){\line(1,-1){10}}
\textcolor{red}{\qbezier(15,0)(22.5,5)(30,10)}
\put(30,10){\line(-1,-2){5}}
\put(30,10){\line(1,-2){5}}
\textcolor{red}{\qbezier(30,10)(37.5,5)(45,0)}
\put(60,10){\line(-1,-2){5}}
\textcolor{red}{\qbezier(60,10)(62.5,5)(65,0)}
\put(40,30){\circle*{2}} \put(40,30){\makebox(0,7)[c]{1}}
\put(30,20){\circle*{2}} \put(30,20){\makebox(0,7)[c]{4}}
\put(40,20){\circle*{2}}\put(40,20){\makebox(0,-7)[c]{8}}
\put(50,20){\circle*{2}}\put(50,20){\makebox(2,7)[c]{3}}
\put(20,10){\circle*{2}}\put(20,10){\makebox(0,7)[c]{9}}
\put(30,10){\circle*{2}}\put(30,10){\makebox(7,2)[c]{14}}
\put(50,10){\circle*{2}}\put(50,10){\makebox(0,-7)[c]{5}}
\put(60,10){\circle*{2}}\put(60,10){\makebox(6,0)[c]{10}}
\put(15,0){\circle*{2}} \put(15,0){\makebox(0,-7)[c]{2}}
\put(25,0){\circle*{2}} \put(25,0){\makebox(0,-7)[c]{11}}
\put(35,0){\circle*{2}} \put(35,0){\makebox(0,-7)[c]{12}}
\put(45,0){\circle*{2}} \put(45,0){\makebox(0,-7)[c]{7}}
\put(55,0){\circle*{2}} \put(55,0){\makebox(0,-7)[c]{13}}
\put(65,0){\circle*{2}} \put(65,0){\makebox(0,-7)[c]{6}}
\put(40,-12){\makebox(0,0)[c]{$T$}}
\put(82.5,10){\makebox(0,0)[c]{\Large{$\Longrightarrow$}}}
\thinlines
\put(105,25){\line(0,1){10}}
\put(110,15){\line(0,-1){10}}
\put(130,25){\line(0,-1){10}}
\put(140,15){\line(0,-1){10}}
\linethickness{0.5mm}
\put(125,35){\line(-1,-2){5}}
\put(125,35){\line(1,-2){5}}
\put(105,25){\line(-1,-2){5}}
\textcolor{red}{\qbezier(105,25)(107.5,20)(110,15)}
\put(130,25){\line(1,-1){10}}
\textcolor{red}{\qbezier(100,5)(105,10)(110,15)}
\put(110,5){\line(-1,-2){5}}
\put(110,5){\line(1,-2){5}}
\textcolor{red}{\qbezier(110,5)(117.5,0)(125,-5)}
\put(140,5){\line(-1,-2){5}}
\textcolor{red}{\qbezier(140,5)(142.5,0)(145,-5)}
\put(105,35){\circle*{2}} \put(105,35){\makebox(6,0)[c]{3}}
\put(125,35){\makebox(0,0)[c]{$\rule{2mm}{2mm}$}}
\put(105,25){\makebox(0,0)[c]{$\rule{2mm}{2mm}$}}
\put(120,25){\circle*{2}}\put(120,25){\makebox(0,-7)[c]{7}}
\put(130,25){\circle*{2}}\put(130,25){\makebox(2,7)[c]{2}}
\put(100,15){\circle*{2}}\put(100,15){\makebox(-2,7)[c]{8}}
\put(110,15){\circle*{2}}\put(110,15){\makebox(7,0)[c]{13}}
\put(130,15){\circle*{2}}\put(130,15){\makebox(0,-7)[c]{4}}
\put(140,15){\circle*{2}}\put(140,15){\makebox(6,0)[c]{9}}
\put(100,5){\circle*{2}} \put(100,5){\makebox(0,-7)[c]{1}}
\put(110,5){\makebox(0,0)[c]{$\rule{2mm}{2mm}$}}
\put(105,-5){\circle*{2}} \put(105,-5){\makebox(0,-7)[c]{10}}
\put(115,-5){\circle*{2}} \put(115,-5){\makebox(0,-7)[c]{11}}
\put(125,-5){\circle*{2}} \put(125,-5){\makebox(0,-7)[c]{6}}
\put(140,5){\makebox(0,0)[c]{$\rule{2mm}{2mm}$}}
\put(135,-5){\circle*{2}} \put(135,-5){\makebox(0,-7)[c]{12}}
\put(145,-5){\circle*{2}} \put(145,-5){\makebox(0,-7)[c]{5}}
\put(122.5,-15){\makebox(0,0)[c]{$\theta(T)$}}
\end{picture}
}
\end{center}
\caption{An example for the bijection $\theta$.\label{fig:half0}}
\end{figure}

Let $\mathcal H_{n,k}$ denote the set of forests of
half-mobile trees on $[n]$ with $k$ improper edges.
The following statement follows immediately
from Theorem \ref{thm:qnkxt-1} and Proposition \ref{prop:theta}.

\begin{thm}For $n\geq 1$, there holds
\begin{align*}
Q_{n,k}(x,t)=\sum_{F\in \mathcal H_{n,k}}x^{\tree(F)-1}t^{\chb(F)}.
\end{align*}
 In other words, we have
\begin{align*}
Q_n(x,y,1,t)=\sum_{F\in\mathcal H_n}x^{\tree(F)-1}
y^{\im(F)}t^{\chb(F)}.
\end{align*}
\end{thm}

\section{Enumeration of plane trees}
We first state a variant of Chu-Vandermonde formula.
\begin{lem}\label{lem:convol}
For $n\geq 1$, there holds
\begin{align}
\sum_{k=0}^{n}{n\choose k}\prod_{i=0}^{k}(x+it)\prod_{j=0}^{n-k-1}(y+jt)
=x\prod_{k=1}^{n}(x+y+kt). \label{eq:convol}
\end{align}
\end{lem}
\pf Write the left-hand side of \eqref{eq:convol} as
$$
(-1)^nn!xt^n\sum_{k=0}^{n}{-x/t-1\choose k}{-y/t\choose n-k}
$$
and then apply the
Chu-Vandermonde convolution formula $\sum_{k=0}^n{a\choose k}{b\choose n-k}={a+b\choose n}$.
\qed

Let $T$ be a plane tree containing the edge $(i,j)$. We define a mapping
 $T\mapsto T'$, called $(i,j)$-\emph{contraction},
by contracting
the edge $(i,j)$ to the vertex $i$ and moving all the children of $j$ to $i$ such that
if $a_1,\ldots, a_r$ (respectively $b_1,\ldots, b_s$) are all the children of $i$ to the
\emph{left} (respectively \emph{right}) of $j$ and $c_1,\ldots, c_t$ are all the children
of $j$, then the children of $i$ in $T'$ are ordered as
$a_1,\ldots, a_r,c_1,\ldots, c_t,b_1,\ldots, b_s.$
An illustration of this contraction is given in Figure~\ref{fig:contract}.

\setlength{\unitlength}{1.2mm}
\begin{figure}[ht]
\begin{center}
{\footnotesize
\begin{picture}(90,33)(5,0)
\put(20,30){\line(0,-1){10}}
\put(20,20){\line(-1,-4){2.5}}
\put(20,20){\line(2,-3){6.6666}}
\put(3,9.3){\makebox(0,0)[c]{$A$}}

\put(-5,8.2){\makebox(0,0)[c]{
\setlength{\unitlength}{2.8mm}
\qbezier(5.8,2.0)(5.8,2.3728)(4.9799,2.6364)
 \qbezier(4.9799,2.6364)(4.1598,2.9)(3.0,2.9)
 \qbezier(3.0,2.9)(1.8402,2.9)(1.0201,2.6364)
 \qbezier(1.0201,2.6364)(0.2,2.3728)(0.2,2.0)
 \qbezier(0.2,2.0)(0.2,1.6272)(1.0201,1.3636)
 \qbezier(1.0201,1.3636)(1.8402,1.1)(3.0,1.1)
 \qbezier(3.0,1.1)(4.1598,1.1)(4.9799,1.3636)
 \qbezier(4.9799,1.3636)(5.8,1.6272)(5.8,2.0) }}

\put(20,20){\line(1,-1){10}}
\put(20,20){\line(-2,-1){20}}
\put(20,20){\line(-4,-3){13.3333}}
\put(32,9){\makebox(0,0)[c]{$B$}}

\put(24,8){\makebox(0,0)[c]{
\setlength{\unitlength}{3.1mm}
\qbezier(5.8,2.0)(5.8,2.3728)(4.9799,2.6364)
 \qbezier(4.9799,2.6364)(4.1598,2.9)(3.0,2.9)
 \qbezier(3.0,2.9)(1.8402,2.9)(1.0201,2.6364)
 \qbezier(1.0201,2.6364)(0.2,2.3728)(0.2,2.0)
 \qbezier(0.2,2.0)(0.2,1.6272)(1.0201,1.3636)
 \qbezier(1.0201,1.3636)(1.8402,1.1)(3.0,1.1)
 \qbezier(3.0,1.1)(4.1598,1.1)(4.9799,1.3636)
 \qbezier(4.9799,1.3636)(5.8,1.6272)(5.8,2.0) }}

\put(17.5,10){\line(-2,-3){6.6666}}
\put(17.5,10){\line(0,-1){10}}
\put(20,20){\line(1,-1){10}}
\put(15,-0.7){\makebox(0,0)[c]{$C$}}

\put(7,-1.7){\makebox(0,0)[c]{
\setlength{\unitlength}{2.5mm}
\qbezier(5.8,2.0)(5.8,2.3728)(4.9799,2.6364)
\qbezier(4.9799,2.6364)(4.1598,2.9)(3.0,2.9)
\qbezier(3.0,2.9)(1.8402,2.9)(1.0201,2.6364)
\qbezier(1.0201,2.6364)(0.2,2.3728)(0.2,2.0)
\qbezier(0.2,2.0)(0.2,1.6272)(1.0201,1.3636)
\qbezier(1.0201,1.3636)(1.8402,1.1)(3.0,1.1)
\qbezier(3.0,1.1)(4.1598,1.1)(4.9799,1.3636)
\qbezier(4.9799,1.3636)(5.8,1.6272)(5.8,2.0) }}

\put(85,30){\line(0,-1){10}}
\put(85,20){\line(2,-3){6.6666}}
\put(68,9.3){\makebox(0,0)[c]{$A$}}

\put(60,8.2){\makebox(0,0)[c]{
\setlength{\unitlength}{2.8mm}
\qbezier(5.8,2.0)(5.8,2.3728)(4.9799,2.6364)
 \qbezier(4.9799,2.6364)(4.1598,2.9)(3.0,2.9)
 \qbezier(3.0,2.9)(1.8402,2.9)(1.0201,2.6364)
 \qbezier(1.0201,2.6364)(0.2,2.3728)(0.2,2.0)
 \qbezier(0.2,2.0)(0.2,1.6272)(1.0201,1.3636)
 \qbezier(1.0201,1.3636)(1.8402,1.1)(3.0,1.1)
 \qbezier(3.0,1.1)(4.1598,1.1)(4.9799,1.3636)
 \qbezier(4.9799,1.3636)(5.8,1.6272)(5.8,2.0) }}

\put(85,20){\line(1,-1){10}}
\put(85,20){\line(-2,-1){20}}
\put(85,20){\line(-4,-3){13.3333}}
\put(97,9){\makebox(0,0)[c]{$B$}}

\put(89,8){\makebox(0,0)[c]{
\setlength{\unitlength}{3.1mm}
\qbezier(5.8,2.0)(5.8,2.3728)(4.9799,2.6364)
 \qbezier(4.9799,2.6364)(4.1598,2.9)(3.0,2.9)
 \qbezier(3.0,2.9)(1.8402,2.9)(1.0201,2.6364)
 \qbezier(1.0201,2.6364)(0.2,2.3728)(0.2,2.0)
 \qbezier(0.2,2.0)(0.2,1.6272)(1.0201,1.3636)
 \qbezier(1.0201,1.3636)(1.8402,1.1)(3.0,1.1)
 \qbezier(3.0,1.1)(4.1598,1.1)(4.9799,1.3636)
 \qbezier(4.9799,1.3636)(5.8,1.6272)(5.8,2.0) }}

\put(85,20){\line(-2,-3){6.6666}}
\put(85,20){\line(0,-1){10}}
\put(82.5,9.3){\makebox(0,0)[c]{$C$}}

\put(74.5,8.3){\makebox(0,0)[c]{
\setlength{\unitlength}{2.5mm}
\qbezier(5.8,2.0)(5.8,2.3728)(4.9799,2.6364)
\qbezier(4.9799,2.6364)(4.1598,2.9)(3.0,2.9)
\qbezier(3.0,2.9)(1.8402,2.9)(1.0201,2.6364)
\qbezier(1.0201,2.6364)(0.2,2.3728)(0.2,2.0)
\qbezier(0.2,2.0)(0.2,1.6272)(1.0201,1.3636)
\qbezier(1.0201,1.3636)(1.8402,1.1)(3.0,1.1)
\qbezier(3.0,1.1)(4.1598,1.1)(4.9799,1.3636)
\qbezier(4.9799,1.3636)(5.8,1.6272)(5.8,2.0) }}

\put(20,20){\circle*{2}}\put(20,20){\makebox(-4,4)[c]{$i$}}
\put(17.5,10){\circle*{2}}\put(17.5,10){\makebox(-4,4)[c]{$j$}}
\put(85,20){\circle*{2}}\put(85,20){\makebox(-4,4)[c]{$i$}}
\put(52.5,15){\makebox(0,0)[c]{{\Large$\Longrightarrow$}}}

\end{picture}
}
\end{center}
\caption{The $(i,j)$-contraction.\label{fig:contract}}
\end{figure}
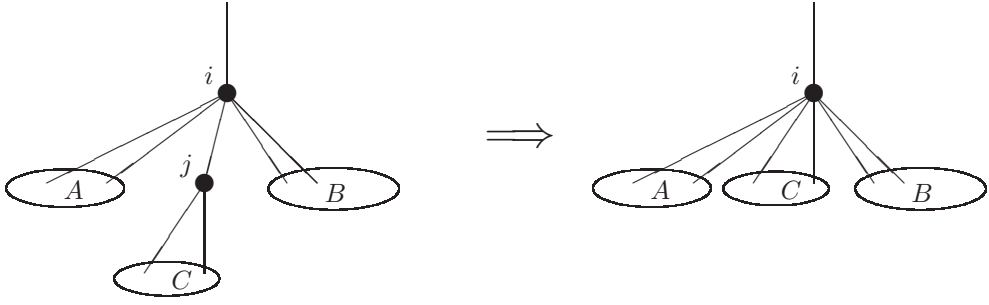

Two plane trees $T_1$ and $T_2$ are said to be {\it $i$-equivalent} and
denoted $T_1\simeq T_2$, if
$T_2$ can be obtained from $T_1$ by reordering the children of the vertex $i$
(see Figure~\ref{fig:2e}).
Moreover, $T_1$ and $T_2$ are said to be {\it $(i,j)$-equivalent} and denoted
$T_1\sim T_2$, if both $T_1$ and $T_2$ contain the edge $(i,j)$ and their images
under the $(i,j)$-contraction are $i$-equivalent.
For instance, the two plane trees $T_3$ and $T_4$ in Figure~\ref{fig:25e}
are $(2,5)$-equivalent,
since after the $(2,5)$-contraction they become respectively $T_1$ and $T_2$
in Figure~\ref{fig:2e}.

\setlength{\unitlength}{1mm}
\begin{figure}[ht]
\begin{center}
{\footnotesize
\begin{picture}(65,39)(0,3)
\put(5,35){\line(-1,-2){5}}
\put(5,35){\line(1,-2){10}}
\put(10,25){\line(-1,-2){5}}
\put(10,25){\line(-0,-1){10}}
\put(10,25){\line(1,-2){5}}
\put(5,15){\line(-1,-2){5}}
\put(5,15){\line(0,-1){10}}
\put(15,15){\line(1,-2){5}}
\put(5,35){\circle*{2}}\put(5,35){\makebox(0,7)[c]{3}}
\put(0,25){\circle*{2}}\put(0,25){\makebox(0,-7)[c]{9}}
\put(10,25){\circle*{2.7}}\put(10,25){\makebox(2,7)[c]{2}}
\put(5,15){\circle*{2}}\put(5,15){\makebox(-6,0)[c]{6}}
\put(10,15){\circle*{2}}\put(10,15){\makebox(0,-7)[c]{8}}
\put(15,15){\circle*{2}}\put(15,15){\makebox(6,0)[c]{1}}
\put(0,5){\circle*{2}}\put(0,5){\makebox(0,-7)[c]{4}}
\put(5,5){\circle*{2}}\put(5,5){\makebox(0,-7)[c]{7}}
\put(20,5){\circle*{2}}\put(20,5){\makebox(0,-7)[c]{10}}
\put(50,35){\line(-1,-2){5}}
\put(50,35){\line(1,-2){10}}
\put(55,25){\line(-1,-2){5}}
\put(55,25){\line(-0,-1){10}}
\put(55,25){\line(1,-2){5}}
\put(50,15){\line(-1,-2){5}}
\put(55,15){\line(0,-1){10}}
\put(55,15){\line(1,-2){5}}
\put(50,35){\circle*{2}}\put(50,35){\makebox(0,7)[c]{3}}
\put(45,25){\circle*{2}}\put(45,25){\makebox(0,-7)[c]{9}}
\put(55,25){\circle*{2.7}}\put(55,25){\makebox(2,7)[c]{2}}
\put(50,15){\circle*{2}}\put(50,15){\makebox(-6,0)[c]{1}}
\put(55,15){\circle*{2}}\put(55,15){\makebox(-4,-7)[c]{6}}
\put(60,15){\circle*{2}}\put(60,15){\makebox(6,0)[c]{8}}
\put(45,5){\circle*{2}}\put(45,5){\makebox(0,-7)[c]{10}}
\put(55,5){\circle*{2}}\put(55,5){\makebox(0,-7)[c]{4}}
\put(60,5){\circle*{2}}\put(60,5){\makebox(0,-7)[c]{7}}
\put(-8,21){\makebox(0,0)[c]{$T_1$}}
\put(68,21){\makebox(0,0)[c]{$T_2$}}
\end{picture}
}
\end{center}
\caption{Two $2$-equivalent plane trees.\label{fig:2e}}
\end{figure}
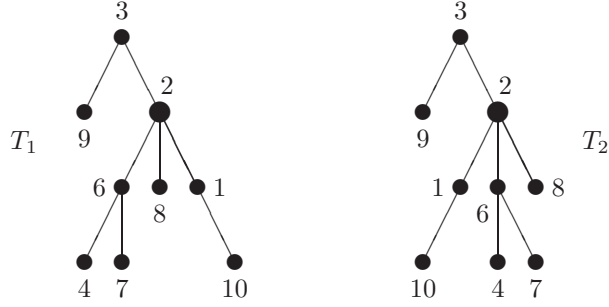

\setlength{\unitlength}{1mm}
\begin{figure}[ht]
\begin{center}
{\footnotesize
\begin{picture}(65,50)(0,-7)
\put(5,35){\line(-1,-2){5}}
\put(5,35){\line(1,-2){10}}
\put(10,25){\line(-1,-2){5}}
\put(10,25){\line(-0,-1){10}}
\put(10,25){\line(1,-2){5}}
\put(5,15){\line(-1,-2){5}}
\put(5,15){\line(0,-1){10}}
\put(10,15){\line(0,-1){10}}
\put(15,15){\line(1,-2){5}}
\put(5,35){\circle*{2}}\put(5,35){\makebox(0,7)[c]{3}}
\put(0,25){\circle*{2}}\put(0,25){\makebox(0,-7)[c]{9}}
\put(10,25){\circle*{2.7}}\put(10,25){\makebox(2,7)[c]{2}}
\put(5,15){\circle*{2}}\put(5,15){\makebox(-6,0)[c]{6}}
\put(10,15){\circle*{2}}\put(10,15){\makebox(4,-7)[c]{5}}
\put(15,15){\circle*{2}}\put(15,15){\makebox(6,0)[c]{1}}
\put(0,5){\circle*{2}}\put(0,5){\makebox(0,-7)[c]{4}}
\put(5,5){\circle*{2}}\put(5,5){\makebox(0,-7)[c]{7}}
\put(10,5){\circle*{2}}\put(10,5){\makebox(0,-7)[c]{8}}
\put(20,5){\circle*{2}}\put(20,5){\makebox(0,-7)[c]{10}}
\put(50,35){\line(-1,-2){5}}
\put(50,35){\line(1,-2){10}}
\put(55,25){\line(-1,-2){10}}
\put(55,25){\line(1,-2){5}}
\put(60,15){\line(-1,-2){10}}
\put(60,15){\line(1,-2){5}}
\put(55,5){\line(1,-2){5}}
\put(50,35){\circle*{2}}\put(50,35){\makebox(0,7)[c]{3}}
\put(45,25){\circle*{2}}\put(45,25){\makebox(0,-7)[c]{9}}
\put(55,25){\circle*{2.7}}\put(55,25){\makebox(2,7)[c]{2}}
\put(50,15){\circle*{2}}\put(50,15){\makebox(-6,0)[c]{1}}
\put(45,5){\circle*{2}}\put(45,5){\makebox(0,-7)[c]{10}}
\put(60,15){\circle*{2}}\put(60,15){\makebox(6,0)[c]{5}}
\put(50,-5){\circle*{2}}\put(50,-5){\makebox(0,-7)[c]{4}}
\put(60,-5){\circle*{2}}\put(60,-5){\makebox(0,-7)[c]{7}}
\put(55,5){\circle*{2}}\put(55,5){\makebox(-6,0)[c]{6}}
\put(65,5){\circle*{2}}\put(65,5){\makebox(0,-7)[c]{8}}
\put(-8,16){\makebox(0,0)[c]{$T_3$}}
\put(73,16){\makebox(0,0)[c]{$T_4$}}
\end{picture}
}
\end{center}
\caption{Two $(2,5)$-equivalent plane trees in $\mathcal P_{10}$.\label{fig:25e}}
\end{figure}
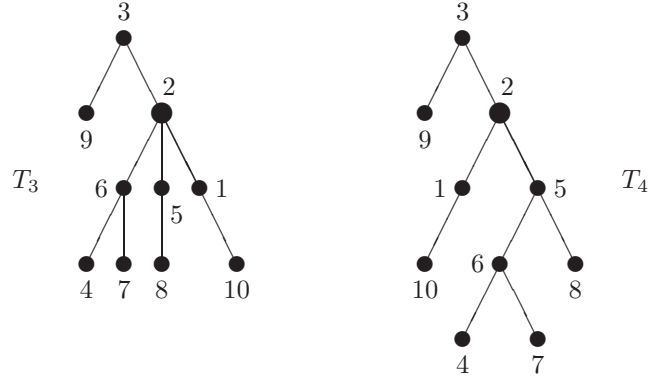

\begin{lem}\label{lem:xixj}
Let $T_0$ be a plane tree on $[n]$ containing the edge $(i, j)$.
Let ${\mathcal T}_0(i,j)$ be the set of plane trees obtained from
the $(i,j)$-contraction of those $i$-equivalent to $T_0$. Then
\begin{align*}
\sum_{T\sim T_0}t^{\eld(T)}\prod_{k=1}^{n}x_k^{\del_T(k)}
=x_i\sum_{T\in{\mathcal T}_0(i,j)}(x_i+x_j+t)^{\del_T(i)}t^{\eld(T)}
\prod_{k\in [n]\setminus\{i,j\}}x_k^{\del_T(k)}.
\end{align*}
\end{lem}

\pf
 Suppose that $\deg_{T_0}(i)+\deg_{T_0}(j)=m+1$. By Eq.~\eqref{eq:general} and
Lemma \ref{lem:convol}, we have
\begin{align*}
\sum_{T\simeq T_0}t^{\eld_T(i)}x_i^{\del_T(i)}t^{\eld_T(j)}x_j^{\del_T(j)}
&=\sum_{k=0}^{m}\sum_{T\sim T_0\atop \deg_T(i)=k+1}t^{\eld_T(i)}x_i^{\del_T(i)}t^{\eld_T(j)}x_j^{\del_T(j)}\\
&=\sum_{k=0}^{m}{m\choose k}\prod_{r=0}^{k}(x_i+rt)
\cdot \prod_{s=0}^{m-k-1}(x_j+st)\\
&=x_i\prod_{k=1}^{m}(x_i+x_j+kt)\\
&=x_i\sum_{T\in{\mathcal T}_0(i,j)}(x_i+x_j+t)^{\del_T(i)}t^{\eld_T(i)}.
\end{align*}
Multiplying by $\prod_{k\in [n]\setminus\{i,j\}}t^{\eld_{T_0}(k)}x_k^{\del_{T_0}(k)}$ we get
the desired result. \qed

Let $\mathcal P_n$ denote the set of all plane trees on $[n]$.
\begin{thm}\label{thm:pplane}
For $n\geq 1$, there holds
\begin{align}
\sum_{T\in \mathcal P_n}t^{\eld(T)}\prod_{i=1}^{n}x_i^{\del_T(i)}
=\prod_{k=0}^{n-2}(x_1+\cdots+x_n+kt).
\label{eq:gen-plane}
\end{align}
\end{thm}

\pf We shall prove \eqref{eq:gen-plane} by induction on $n$. The identity
obviously holds for $n=1$. Suppose \eqref{eq:gen-plane} holds for $n-1$.
For $ i,j\in [ n]$ ($i\neq j$),
let $\mathcal P_{n}(i,j)$ be the set of all plane trees on $[n]$ containing
the edge $(i,j)$. Let $\mathcal A$ be a maximal set of plane trees in $\mathcal P_{n}(i,j)$
that are pairwise not $(i,j)$-equivalent and $\mathcal B$ the set of plane trees
on $[n]\setminus\{j\}$ obtained from the $(i,j)$-contraction of those in $\mathcal A$.

By Lemma \ref{lem:xixj}, we have
\begin{align*}
\sum_{T\in \mathcal P_{n}(i,j)}t^{\eld(T)} \prod_{k=1}^{n}x_k^{\del_T(k)}
&=\sum_{T_0\in \mathcal A}
\sum_{T\sim T_0}t^{\eld(T)}\prod_{k=1}^{n}x_k^{\del_T(k)} \\
&=x_i\sum_{T_0'\in \mathcal B} \sum_{T\simeq T_0'}
(x_i+x_j+t)^{\del_T(i)}t^{\eld(T)}
\prod_{k\in [n]\setminus\{i,j\}}x_k^{\del_T(k)}.
\end{align*}
By the induction hypothesis, the last double sum is equal to
\begin{align*}
&x_i\sum_{T\in\mathcal P_{[n]\setminus\{j\}}}
(x_i+x_j+t)^{\del_T(i)}t^{\eld(T)}
\prod_{k\in [n]\setminus\{i,j\}}x_k^{\del_T(k)} \\
&=x_i\prod_{k=0}^{n-3}(x_1+\cdots+x_n+(k+1)t).
\end{align*}
Summing over all of the pairs $(i,j)$ with $i\neq j$ we get
\begin{align}
&\hskip -3mm \sum_{\substack{1\leq i,j\leq n\\ i\neq j}}
\sum_{T\in \mathcal P_{n}(i,j)}t^{\eld(T)} \prod_{k=1}^{n}x_k^{\del_T(k)} \nonumber\\
&=(n-1)(x_1+\cdots+x_n)\prod_{k=0}^{n-3}(x_1+\cdots+x_n+(k+1)t) \nonumber\\
&=(n-1)\prod_{k=0}^{n-2}(x_1+\cdots+x_n+kt).\label{eq:pnij1}
\end{align}
Noticing that any plane tree on $[n]$ has $n-1$ edges, we have
\begin{equation}\label{eq:pnij2}
\sum_{\substack{1\leq i,j\leq n\\ i\neq j}}
\sum_{T\in \mathcal P_{n}(i,j)}t^{\eld(T)}\prod_{k=1}^{n}x_k^{\del_T(k)}
=(n-1)\sum_{T\in \mathcal P_n}t^{\eld(T)}\prod_{k=1}^{n}x_k^{\del_T(k)}.
\end{equation}
Comparing \eqref{eq:pnij1} and \eqref{eq:pnij2} yields the desired
formula for $\mathcal P_n$.
\qed

\begin{cor}The number of unlabeled plane trees on $n+1$ vertices
is the Catalan number $C_n=\frac{1}{n+1}{2n\choose n}$.
\end{cor}
\pf Setting $x_1=\cdots=x_n=t=1$ in \eqref{eq:gen-plane}
gives the number of labeled plane trees
on $n$ vertices, i.e., $|\mathcal P_n|=(2n-2)!/(n-1)!$.
We then obtain the number of unlabeled plane trees on $n$ vertices by
dividing
$|\mathcal P_{n}|$ by $n!$.
\qed

\begin{cor}The number of unlabeled plane trees with $k$ leaves on
$n+1$ vertices is the Narayana number $N_{n,k}=\frac{1}{n}{n\choose k}{n\choose k-1}$.
\end{cor}
\pf Setting $t=1$ in \eqref{eq:gen-plane} and replacing $n$ by $n+1$, we have
\begin{align*}
\sum_{T\in \mathcal P_{n+1}}\prod_{i=1}^{n+1}x_i^{\del_T(i)}
=\prod_{s=0}^{n-1}(x_1+\cdots+x_{n+1}+s).
\end{align*}
Note that the vertex $i$ is a leaf of the plane tree $T$ if and only if $\del_T(i)=0$.
Hence the number $a_{n,k}$ of plane trees on $[n+1]$ with leaves $1,2,\ldots,k$ is
the sum of the coefficients of monomials $x_{k+1}^{r_{k+1}}\cdots x_{n+1}^{r_{n+1}}$
with $r_j\geq 1$ for all $k+1\leq j\leq n+1$ in the expansion of the symmetric polynomial
$$
P=\prod_{s=0}^{n-1}(x_{k+1}+\cdots+x_{n+1}+s).
$$
Note that the sum of the coefficients of monomials in $P$ not containing
 $i$ ($0\leq i\leq n+1-k$) fixed variables is
$\prod_{s=0}^{n-1}(n+1-k-i+s)$. By the Principle of Inclusion-Exclusion and
the Chu-Vandermonde convolution formula,
one sees that
\begin{align*}
a_{n,k}&=\sum_{i=0}^{n-k+1}(-1)^i{n-k+1\choose i}\prod_{s=0}^{n-1}(n+1-k-i+s)\\
&=n!\sum_{i=0}^{n-k+1}(-1)^i{n-k+1\choose i}{2n-k-i\choose n}\\
&=n!{n-1\choose k-1}.
\end{align*}
 Therefore,
the number of unlabeled plane trees with $k$ leaves on
$n+1$ vertices is equal to
\begin{align*}
\frac{a_{n,k}}{k!(n-k+1)!}=\frac{n!}{k!(n-k+1)!}{n-1\choose k-1}=N_{n,k}.
\tag*{\qed}
\end{align*}

Let $\mathcal P_{n}^{(r)}$ denote the set of
plane trees on $[n]$ with a specific root $r\in [n]$.
Then we can refine Theorem~\ref{thm:pplane} as follows.

\begin{thm}\label{thm:refine} For $n\geq 1$, there holds
\begin{align}
\sum_{T\in \mathcal P_{n}^{(r)}}t^{\eld(T)}\prod_{i=1}^{n}x_i^{\del_T(i)}
=x_r\prod_{k=1}^{n-2}(x_1+\cdots+x_{n}+kt). \label{eq:gen-on}
\end{align}
\end{thm}

\pf Suppose Eq.~\eqref{eq:gen-on} hold for $n-1$.
Let $\mathcal P_n^{(r)}(i,j)$ denote the subset of plane trees
in $\mathcal P_n(i,j)$ with root $r$. Similarly to the proof of
Theorem~\ref{thm:pplane}, we can show that
\begin{align*}
\sum_{T\in \mathcal P_n^{(r)}(i,j)}t^{\eld(T)}\prod_{i=1}^{n}x_i^{\del_T(i)}
=\begin{cases}
x_rx_i\prod_{k=1}^{n-3}(x_1+\cdots+x_{n}+(k+1)t), &\text{if $i\neq r$},\\
(x_r+x_j+t)x_i\prod_{k=1}^{n-3}(x_1+\cdots+x_{n}+(k+1)t), &\text{if $i=r$.}
\end{cases}
\end{align*}
The proof then follows from computing the following sum:
\begin{align*}
\sum_{\substack{1\leq i,j\leq n\\ i\neq j}}
\sum_{T\in \mathcal P_n^{(r)}(i,j)}t^{\eld(T)}\prod_{i=1}^{n}x_i^{\del_T(i)},
\end{align*}
and using the fact that any plane trees on $[n]$ has $n-1$ edges. \qed

\begin{cor}\label{cor:rsroot}
For $1\leq r<s\leq n$, there holds
\begin{align*}
x_r^{-1}\sum_{T\in \mathcal P_{n}^{(r)}}t^{\eld(T)}\prod_{i=1}^{n}x_i^{\del_T(i)}
=x_s^{-1}\sum_{T\in \mathcal P_{n}^{(s)}}t^{\eld(T)}\prod_{i=1}^{n}x_i^{\del_T(i)}.
\end{align*}
\end{cor}

It would be interesting to have a combinatorial proof of the above identity.
We give such a proof in the case $(r,s)=(1,2)$.
Let $T$ be any plane tree in $\mathcal P_{n}^{(1)}$.
Suppose all the children of the root $1$ are
$a_1,\ldots, a_m$ and $2$ is a descendant of $a_t$.
Hence, $a_1,\ldots, a_{t-1}$ are elder, while $a_t$ is younger.
Assume all the children of $2$ are $b_1,\ldots,b_l$.
Exchanging the subtrees with roots $a_{t+1},\ldots,a_m$ in $T$ and the
subtrees with roots $b_1,\ldots,b_l$ in their previous orders,
and then exchanging the labels $1$ and $2$, we obtain a plane tree
$T'\in \mathcal P_{n}^{(2)}$.
It is easy to see that the mapping $T\mapsto T'$ is a bijection from
$\mathcal P_n^{(1)}$ to $\mathcal P_n^{(2)}$. Moreover, one sees that
$\eld(T)=\eld(T')$, $\del_T(1)=\del_{T'}(1)+1$, $\del_T(2)=\del_{T'}(2)-1$,
and $\del_T(i)=\del_{T'}(i)$ if $i\neq 1,2$.

As applications of Theorem~\ref{thm:refine} we derive three classical results on planted and
plane forests, which correspond, respectively, to Theorem 5.3.4, Corollary 5.3.5 and
Theorem 5.3.10 in Stanley's book~\cite{Stanley99}.

Recall that a {\it planted forest} $\sigma$ (or, {\it rooted forest})
is a graph whose connected components are rooted trees.
If the vertex set of $\sigma$ is $[n]$, then we
define the {\it ordered degree sequence}
$\delta(\sigma)=(d_1,\ldots,d_n)$, where $d_i=\deg(i)$.

\begin{cor}\label{cor:planted}
Let ${\bf d}=(d_1,\ldots,d_n)\in\mathbb{N}^n$ with $\sum d_i=n-k$.
Then the number of planted forests on $[n]$
(necessarily with $k$ components) with ordered degree sequence
${\bf d}$ is given by
$$
{n-1\choose k-1}{n-k\choose d_1,\ldots,d_n}.
$$
\end{cor}

\pf
Equating the coefficients of $x_1^k x_2^{d_1}\cdots x_{n+1}^{d_n}$ in
\eqref{eq:gen-on} with $r=1$, $t=0$ and $n$ replaced by $n+1$,
we get the desired result.
\qed

Given a planted forest $\sigma$, define the {\it type} of $\sigma$
to be the sequence
$$
{\rm type}\,\sigma=(r_0,r_1,\ldots),
$$
where $r_i$ vertices of $\sigma$ have degree $i$. Similarly we can define
the {\it type} for a plane forest.
As is well-known (see \cite[p.~30]{Stanley99}),
Corollary~\ref{cor:planted} can be restated in the following equivalent form.

\begin{cor} \label{cor:type-planted}
Let ${\bf r}=(r_0,\ldots,r_m)\in\mathbb{N}^{m+1}$ with $\sum r_i=n$
and $\sum (1-i)r_i=k>0$.
Then the number of planted forests on $[n]$
(necessarily with $k$ components) of type ${\bf r}$ is given by
$$
{n-1\choose k-1}\frac{(n-k)!}{0!^{r_0}\cdots m!^{r_m}}
{n\choose r_0,\ldots,r_m}.
$$
\end{cor}

 A {\it plane forest} is a family of plane trees
in which the roots of the plane trees are linearly ordered.

\begin{cor} \label{cor:plane}
Let ${\bf r}=(r_0,\ldots,r_m)\in\mathbb{N}^{m+1}$ with $\sum r_i=n$
and $\sum (1-i)r_i=k>0$.
Then the number of unlabeled plane forests on $[n]$
(necessarily with $k$ components) of type ${\bf r}$ is given by
$$
\frac{k}{n}{n\choose r_0,\ldots,r_m}.
$$
\end{cor}

\pf
Let us see how many plane forests we can obtain from
a planted forest $\sigma$ of type ${\bf r}$ by ordering the vertices.
For any vertex $v$ of $\sigma$, there are
$\deg(v)!$ ways to linearly order its children.
By definition, there are $r_i$ vertices of $\sigma$ having degree $i$.
Besides, there are $k!$ ways to linearly order the trees of $\sigma$.
So, we can obtain $0!^{r_0}\cdots m!^{r_m}k!$ plane forests from $\sigma$.
Hence it follows from Corollary~\ref{cor:type-planted} that
the total number of plane forests on $[n]$ of type ${\bf r}$
is equal to
\begin{align}
{n-1\choose k-1}\frac{(n-k)!}{0!^{r_0}\cdots m!^{r_m}}
{n\choose r_0,\ldots,r_m} 0!^{r_0}\cdots m!^{r_m}k!
=(n-1)!k{n\choose r_0,\ldots,r_m}.
\label{eq:ff1}
\end{align}
On the other hand, every unlabeled plane forest with $n$ vertices has
$n!$ different labeling methods. Dividing \eqref{eq:ff1} by $n!$ yields
the desired formula. \qed

\section{Forests of plane trees}
When $x=r$ is an integer, we can give another interpretation
for the polynomial $rQ_{n-r,k}(r,t)$ in terms of forests.
Let $\mathcal F_{n,k}^r$ denote the set of
forests of $r$ plane trees on $[n]$ with $k$ improper edges
and with roots $1,\ldots,r$. The following is a generalization
of a theorem of Shor~\cite{Shor}, which corresponds to the case $t=0$.

\begin{thm}\label{thm-gen-shor}
The generating function for forests of $r$ plane trees on $[n]$ with
$k$ improper edges and with roots $1,\ldots,r$ by number of elder vertices
is $rQ_{n-r,k}(r,t)$. Namely,
\begin{align}
rQ_{n-r,k}(r,t)=\sum_{F\in\mathcal F_{n,k}^r}t^{\eld(F)}.
\label{eq:rqnr}
\end{align}
\end{thm}

\pf We proceed by induction on $n$. Identity \eqref{eq:rqnr} obviously holds
for $n=r+1$. Suppose \eqref{eq:rqnr} holds for $n$.
Similarly to the proof of \ref{thm:qnkxt-1},
we can show that
\begin{align}
\sum_{T\in \mathcal F_{n+1,k}^r[\deg(n+1)=0]}t^{\eld(T)}
&=[n+t(n-r)]\sum_{F\in \mathcal F_{n,k}^r}t^{\eld(T)},
 \label{eq:fdeg=0} \\
\sum_{F\in \mathcal F_{n+1,k}[\deg(n+1)>0]}t^{\eld(T)}
&=kt\sum_{F\in \mathcal F_{n,k}^r}t^{\eld(T)}
 +(n+k-r-1)\sum_{F\in \mathcal P_{n,k-1}^r}t^{\eld(T)}.
 \label{eq:fdeg>0}
\end{align}
Summing \eqref{eq:fdeg=0} and \eqref{eq:fdeg>0} and using the
induction hypothesis and \eqref{rec-qnkxt}, we obtain
\begin{align*}
\sum_{T\in \mathcal F_{n+1,k}^r}t^{\eld(T)}
&=[n+t(n+k-r)]rQ_{n-r,k}(r,t)+(n+k-r-1)rQ_{n-r,k-1}(r,t), \nonumber\\
&=rQ_{n+1-r,k}(r,t).
\end{align*}
Namely, Eq.~\eqref{eq:rqnr} holds for $n+1$. This completes the proof.
\qed

\section{Proof of the duality formula for $Q_n$}
First we restate the duality formula~\eqref{eq:conj} in terms of $Q_{n,k}(x,t)$.
\begin{lem}
For $n\geq 2$, the duality formula~\eqref{eq:conj} is equivalent to
\begin{equation}\label{eq:rec2}
Q_{n,k}(x,t)=(x-k+t+1)Q_{n-1,k}(x+t+1,t)+(n+k-2)Q_{n-1,k-1}(x+t+1,t).
\end{equation}
\end{lem}
\pf
Plugging \eqref{eq:expansion} into \eqref{eq:conj} we see that \eqref{eq:conj} is equivalent to
the following recurrence relation for $Q_{n,k}(x,t)$:
\begin{align}\label{eq:mainconj}
Q_{n,k}(x,t)=Q_{n,k}(-(x+n+nt)/t,1/t)(-t)^{n-k-1}.
\end{align}
Setting $X_n=-(x+n+nt)/t$ and $T=1/t$, by means of \eqref{rec-qnkxt}
 we have
\begin{align}
Q_{n,k}(X_n,T)=-\frac{x-k+t+1}{t}Q_{n-1,k}\left(X_n, T\right)
+(n+k-2)Q_{n-1,k-1}\left(X_n, T\right). \label{eq:new}
\end{align}
Multiplying \eqref{eq:new} by $(-t)^{n-k-1}$,
we see that \eqref{eq:mainconj} is equivalent to \eqref{eq:rec2}. \qed

Now, by replacing $x$ with $x-t-1$ in \eqref{eq:rec2}, we get
\begin{equation}\label{eq:rec3}
Q_{n,k}(x-t-1,t)=(x-k)Q_{n-1,k}(x,t)+(n+k-2)Q_{n-1,k-1}(x,t).
\end{equation}
Subtracting \eqref{eq:rec3} from \eqref{rec-qnkxt}, we are led
to the following equivalent identity:
\begin{lem}
For $n\geq 2$ and $0\leq k\leq n-1$ there holds
\begin{align}
Q_{n,k}(x,t)-Q_{n,k}(x-t-1,t)=(t+1)(n+k-1)Q_{n-1,k}(x,t).
\label{eq:qnkxt-diff}
\end{align}
\end{lem}
\pf We proceed by induction on $n$. Eq.~\eqref{eq:qnkxt-diff} is obviously true
for $n=2$. Suppose it is true for $n-1$. By the definition \eqref{rec-qnkxt}
of $Q_{n,k}(x,t)$ we have
\begin{align*}
&\hskip -2mm Q_{n,k}(x,t)-Q_{n,k}(x-t-1,t) \\
&=[x+n-1+t(n+k-2)]Q_{n-1,k}(x,t)+(n+k-2)Q_{n-1,k-1}(x,t) \\
&\quad -[x+n-2+t(n+k-3)]Q_{n-1,k}(x-t-1,t)-(n+k-2)Q_{n-1,k-1}(x-t-1,t)\\
&=[x+n-2+t(n+k-3)][Q_{n-1,k}(x,t)-Q_{n-1,k}(x-t-1,t)]\\
&\quad +(t+1)Q_{n-1,k}(x,t)+(n+k-2)[Q_{n-1,k-1}(x,t)-Q_{n-1,k-1}(x-t-1,t)],
\end{align*}
and the induction hypothesis implies that the above quantity is equal to
\begin{align*}
&(t+1)[x+n-2+t(n+k-3)](n+k-2)Q_{n-2,k}(x,t)+(t+1)Q_{n-1,k}(x,t) \\
&\quad +(t+1)(n+k-2)(n+k-3)Q_{n-2,k-1}(x,t) \\
&=(t+1)(n+k-2)\{[x+n-2+t(n+k-3)]Q_{n-2,k}(x,t)+(n+k-3)Q_{n-2,k-1}(x,t)\} \\
&\quad +(t+1)Q_{n-1,k}(x,t) \\
&=(t+1)(n+k-2)Q_{n-1,k}(x,t)+(t+1)Q_{n-1,k}(x,t) \\
&=(t+1)(n+k-1)Q_{n-1,k}(x,t).
\end{align*}
Thus \eqref{eq:qnkxt-diff} is true for $n$. This completes the proof. \qed

\begin{rmk}
We can also prove Theorem \ref{thm:conj} by arguing directly with $Q_n$ instead of
$Q_{n,k}$. Indeed, the theorem is equivalent to saying that
\begin{equation}\label{eq:equivthm}
F_{n+1}(x-z-t)=[x+nz+(y-z)(n+y\partial_y)]F_n(x),
\end{equation}
where $F_n(x)=Q_n(x,y,z,t)$.
Subtracting \eqref{eq:equivthm} from the definition~\eqref{eq:recdef} yields
\begin{equation}\label{eq:equivbis}
F_{n+1}(x)-F_{n+1}(x-z-t)=(z+t)(n+y\partial_y)F_n(x).
\end{equation}
By \eqref{eq:recdef} and the induction hypothesis,
the left-hand side of \eqref{eq:equivbis} may be written as
\begin{align*}
&[x+nz+(y+t)(n+y\partial_y)]F_n(x)
 -[x-z-t+nz+(y+t)(n+y\partial_y)]F_n(x-z-t)\\
&=[x+(n-1)z-t][F_n(x)-F_n(x-z-t)]+(z+t)F_n(x)\\
&\quad{} +(y+t)(n+y\partial_y)[F_n(x)-F_n(x-z-t)]\\
&= (z+t)\left\{[x+(n-1)z-t](n-1+y\partial_y)
 +(y+t)(n+y\partial_y)(n-1+y\partial_y)\right\}F_{n-1}(x)\\
&\hspace{3cm}+(z+t)F_n(x).
\end{align*}
Applying the differential operator identity
$$
(y+t)(n+y\partial_y)(n-1+y\partial_y)
=(n-1+y\partial_y)[(y+t)(n+y\partial_y)-y],
$$
to the last expression, we obtain
\begin{align*}
&(z+t)(n-1+y\partial_y)\left\{[x+(n-1)z+ (y+t)(n-1+y\partial_y)]F_{n-1}(x)\right\}
  +(z+t)F_n(x)\\
&=(z+t)(n+y\partial_y)F_n(x),
\end{align*}
as desired.

\end{rmk}

\section{Open problems}
Although we have shown  the fecundity of polynomials $Q_n$ in
 the enumeration of plane trees and
forests, there are still further interesting problems.

First of all,
is there any connection between these polynomials and the Lambert $W$
function or a generalization of the Lambert $W$ function? In particular,
is there an analogue of the formula \eqref{eq:rama} or \eqref{eq:prefun} for $Q_n$?

Secondly, it would be interesting to have a better combinatorial
understanding of the polynomials $Q_n$. For example,
from Theorems \ref{thm:qnkxt-1} and \ref{thm:x-t-1}
we deduce that
\begin{align}
\sum_{T\in \mathcal O_{n+1,k}}x^{\del_T(1)-1}t^{\eld(T)}
=\sum_{T\in \mathcal P_{n,k}}(x+t+1)^{\del_T(1)}t^{\eld(T)}.
\label{eq:equiv}
\end{align}
Is there a direct combinatorial proof of \eqref{eq:equiv}?
Also, the duality formula \eqref{eq:conj} deserves a combinatorial proof.
For $t=0$, such proofs have been given by Chen and Guo~\cite{CG}.

Thirdly, since our proof of Theorems~\ref{thm:pplane} and \ref{thm:refine}
is by induction,  it remains to find direct bijective
proofs of them.

Finally, using a dual statistic, the number of ``improper vertices,"
Gessel and Seo~\cite{GS} have recently given several combinatorial
interpretations of the polynomials
\begin{equation}\label{eq:gs}
xQ_n(x,z,z,t-z)=x\prod_{k=1}^{n-1}(x+(n-k)z+kt).
\end{equation}
On the other hand, we derive from Eqs. \eqref{eq:qnxt}, \eqref{eq:expansion}
and Theorem~\ref{thm:qnkxt-1} that the same polynomials~\eqref{eq:gs}
also have the following expression:
$$
\sum_{T\in\mathcal P_{n+1}^{(1)}}x^{\del_T(1)}(t-z)^{\eld(T)}z^{n-\del_T(1)-\eld(T)}.
$$
It might be of interest to have  a combinatorial
explanation for these different interpretations.

\indent{\bf Acknowledgment.}
We would like to thank our colleague Fr\'ed\'eric Chapoton for
his patience and helpful discussions.

\renewcommand{\baselinestretch}{1}

\end{document}